\documentclass[3p, times, 11pt, numafflabel]{elsarticle}
\usepackage[utf8]{inputenc}
\usepackage[T1]{fontenc}
\usepackage[indent, skip=0.12\baselineskip plus 2pt]{parskip}
\usepackage{setspace}
\biboptions{sort}

\journal{Lect. Notes Comput. Sci. Eng.}

\usepackage{amsmath, amssymb,  babel, booktabs, enumitem, microtype}
\usepackage{algorithm, algpseudocode, tabularx}
\usepackage[group-digits=integer, group-minimum-digits=5]{siunitx}

\newdefinition{definition}{Definition}
\newdefinition{remark}{Definition}

\usepackage{mathtools}
\DeclarePairedDelimiter\abs{\lvert}{\rvert}

\DeclareMathOperator{\curl}{curl}
\DeclareMathOperator{\grad}{grad}
\DeclareMathOperator{\p}{\partial}
\DeclareMathOperator{\R}{\mathbb{R}}
\DeclareMathOperator{\tr}{tr}
\DeclareMathOperator{\vol}{vol}
\renewcommand{\vec}[1]{\boldsymbol{#1}}
\providecommand{\D}{\mathrm{d}}
\providecommand{\E}{\mathrm{e}}
\newcommand{\tetgen}{TetGen}

\usepackage{subcaption}
\DeclareCaptionLabelFormat{boldparens}{\textbf{(#2)}}
\usepackage[colorlinks=true, allcolors=blue, hypertexnames=false]{hyperref}
\usepackage[capitalize]{cleveref} 
\crefname{equation}{}{}
\Crefname{equation}{Equation}{Equations}
\crefname{figure}{Fig.}{Figs.}
\Crefname{figure}{Figure}{Figures}
\crefname{section}{Sect.}{Sects.}
\Crefname{section}{Section}{Sections}
\crefname{algorithm}{Algorithm}{Algorithms}
\crefname{table}{Table}{Tables}

\newcommand{\runinhead}[1]{\subsection*{#1}}

\begin{document}

\begin{frontmatter}

\date{\today}

\title{Adaptive Mesh Refinement for~Electromagnetic Simulation}

\author[CCAS]{Alexey Belokrys-Fedotov}
\ead{belokrys.fedotov@yandex.ru}

\author[CCAS]{Vladimir Garanzha}
\ead{garan@ccas.ru}%

\author{Lennard Kamenski}
\ead{l.kamenski@arcor.de}

\author[HW]{Alexandr Chikitkin}
\ead{chikitkin.aleksandr1@huawei.com}

\author[HW]{Evgeniy Pesnya}
\ead{pesnya.evgeniy@huawei-partners.com}

\author[HW]{Nikita Aseev}
\ead{aseev.nikita1@huawei.com}

\author[HW]{Andrey Vorobyev}
\ead{andrey.vorobyev@huawei.com}%

\address[CCAS]{Dorodnicyn Computing Center FRC CSC RAS}
\address[HW]{Huawei Technologies}

\begin{abstract}
We consider problems related to initial meshing and adaptive mesh refinement for the 
electromagnetic simulation of various structures.
The quality of the initial mesh and the performance of the adaptive refinement are of great importance for the finite element solution of the Maxwell equations, since they directly affect the accuracy and the computational time.
In this paper, we describe the complete meshing workflow, which allows the simulation of arbitrary structures.
Test simulations confirm that the presented approach allows to reach the quality of the industrial simulation software.
\end{abstract}

\begin{keyword}
Maxwell equations, adaptive mesh refinement, tetrahedral meshes
\MSC[2020] 65N50 
\\[0.6\baselineskip]%
{\footnotesize{%
This is a preprint of the following chapter:
  A.~Belokrys-Fedotov et al.,
  Adaptive mesh refinement for electromagnetic simulation,
published in
  Numerical Geometry, Grid Generation and Scientific Computing, NUMGRID 2022,
edited by
  V.~Garanzha and L.~Kamenski,
2024, Springer, Cham.
Reproduced with permission of Springer Nature Switzerland.
\\%
The final authenticated version is available online at: \url{https://doi.org/10.1007/978-3-031-59652-0_7}.
}}
\end{keyword}

\end{frontmatter}

\section{Introduction}
Numerical simulation of electromagnetic waves with the Finite Element Method (FEM) is a very important tool for the design of modern devices.
Despite the rapid growth of computing power and available memory resources, simulations of real devices are often prohibitively expensive.
Both computational time and memory should be reduced as much as possible, while keeping the accuracy at a reasonable level.
This leads to quite strict requirements for FEM meshes: they should be close to optimal in terms of the number of elements.
Such optimality can only be achieved by adaptive mesh refinement, when the information about the numerical solution on the current mesh is incorporated into the mesh generation and adaptation processes.

Usually, this information is extracted in the form of error estimates, scalar quantities that majorize the exact error or, at least, provide some qualitative information about the true error distribution.
Given error estimates, the adaptive algorithm should refine the current mesh in some optimal way, but at the same time with low computational cost.
These controversial requirements reject many well-known mesh refinement techniques with beautiful theoretical properties, such as bisection and Delaunay refinement, to name a few.
The only way is to combine some parts of classical algorithms with good heuristics.

Another challenge for mesh refinement is how to handle the high aspect ratio and geometric near-degeneracies, which are inherent in industrial applications.

In this paper we address all of the aforementioned aspects of adaptive meshing.
\Cref{sec:test_structures} presents two structures used for testing and highlights their properties which have an impact on the meshing process.
\Cref{sec:equations,sec:fem} briefly introduce the time-harmonic Maxwell equations and FEM for their solution.
\Cref{sec:thin_layers_meshing} describes in detail the key parts of the meshing algorithm: important modifications of the classical Delaunay refinement algorithm, special procedures for sliver removal, and variational vertex smoothing.
The adaptive mesh refinement algorithm is presented in \cref{sec:amr}.
Simulation results obtained with the presented meshing workflow are given in \cref{sec:simulation_results}.

\section{Test Structures}%
\label{sec:test_structures}

Before introducing the meshing technique itself, let us briefly describe two relatively simple but quite representative test cases exhibiting a number of problems that should be resolved by computational technology in order to be applicable to industrial electromagnetic (EM) simulations.

\begin{figure}[t]%
  \centering{}%
  \includegraphics[width=0.8\textwidth]{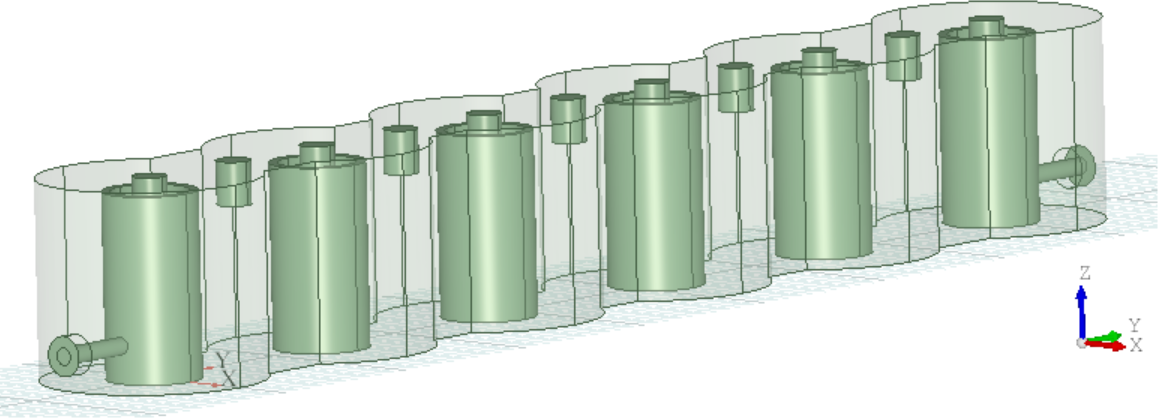}%
  \\[1ex]%
  \includegraphics[width=0.8\textwidth]{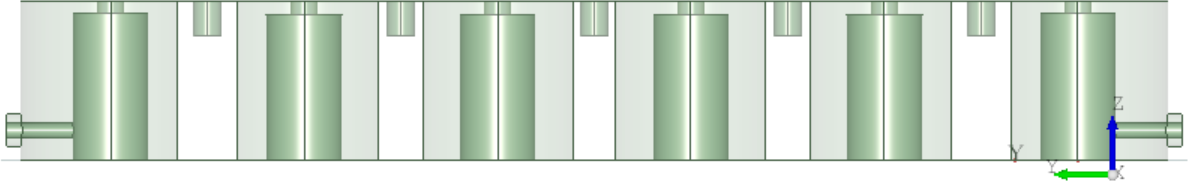}%
  \caption{Cavity resonator filter (6cav)}%
  \label{fig:6cav_structure}%
\end{figure}

\runinhead{Resonator filter with 6 cavities (6cav)}
The first structure is a resonator cavity filter (\cref{fig:6cav_structure}). The interior of this filter is a vacuum, the walls are perfect electrical conductor (PEC) boundaries, and two symmetrically located ports are coaxial waveguide ports.
Mesh adaptation was performed at the frequency of \qty{2.14}{\GHz}.

At first glance, this structure is very simple, but computations show that the accuracy in S-parameters is very sensitive to the mesh density near the cylindrical screws entering the cavities from above. Furthermore, the large area of the boundary makes the point generation rules in the adaptive mesh refinement (AMR) particularly important. Improper point generation can lead to an over-refined surface mesh and a too coarse volume mesh, which, again, is critical for accuracy.

\begin{figure}[t]%
   \includegraphics[width=0.49\textwidth]{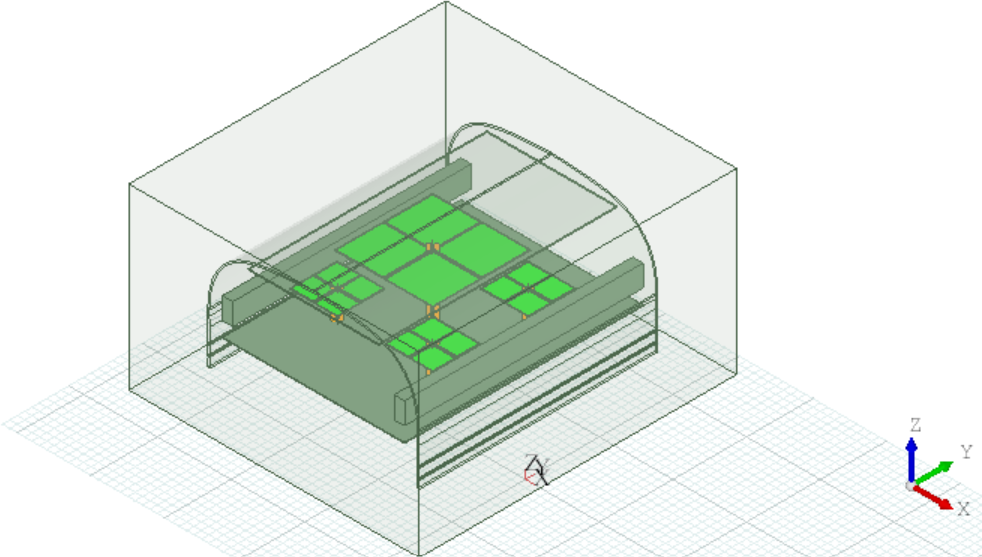}%
   \hfill{}
    \includegraphics[width=0.49\textwidth]{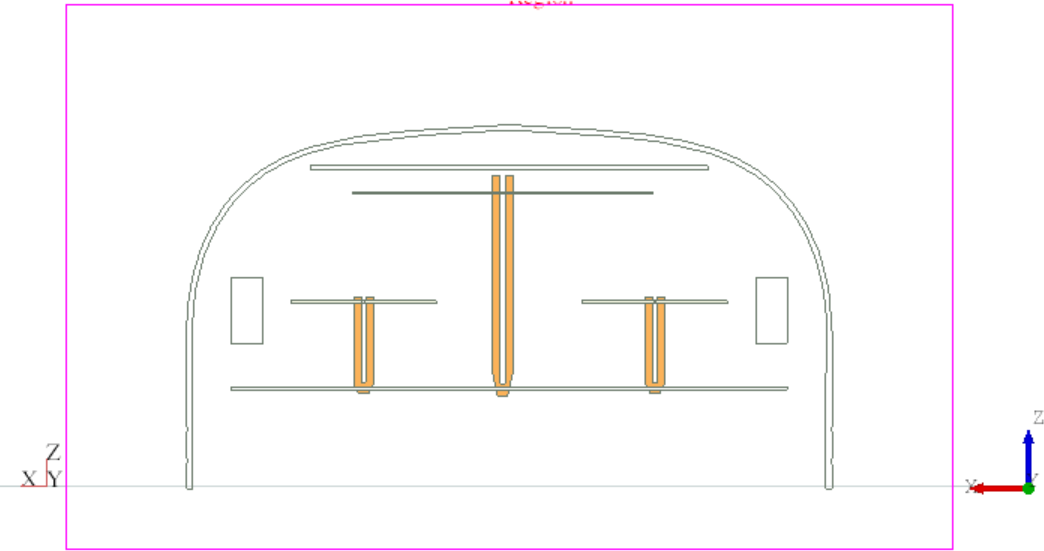}%
  \caption{Antenna array unit (AAU)}%
  \label{fig:aau_struct}%
\end{figure}

\runinhead{Antenna array unit (AAU)}
The second test structure is a small antenna array unit (AAU) (\cref{fig:aau_struct}).
It is placed in a vacuum box that is truncated using first order absorbing boundary conditions \cref{eq:ABC}.
The structure contains a large electrically transparent radome and a dielectric layer with $\varepsilon_r = 2.97$ below the radome. There are four high frequency and one low frequency radiating elements, ground below and two PEC boxes on the sides. Each radiating element consists of two pairs of PEC frames mounted on the dielectric at the top and bottom. For the high frequency elements, $\varepsilon_r = 4.1$ and $\mu_r = 1$. For the low frequency element, $\varepsilon_r = 2.7$ and $\mu_r = 1$. Radiating elements are connected to the ground with flat PEC pins.
Mesh adaptation was performed at the frequency of \qty{2.2}{\GHz}.

The AAU has several peculiarities that strongly affect the meshing process:
\begin{itemize}
  \item Multiple materials, meaning that the surface of the model is not manifold.
  \item A large number of internal surfaces, including the complex case of surfaces inside a single material.
  \item Many thin layers with aspect ratio from $10$ to $100$ (dielectric squares, radome).
  \item Small curvilinear grooves on the radome (left bottom part of the figure).
  \item Geometric flaws resulting from imperfect CAD processing. Namely, in the middle of the radome there is a tiny step between two halves. Such defects can emerge when the whole structure is created by adding an imperfect reflection of one half. Of course, it can be fixed manually, but robust industrial algorithms must handle such cases automatically.
\end{itemize}

\section{Governing Equations and Electromagnetic Simulation Cases}%
\label{sec:equations}

We consider the source-free time-harmonic Maxwell equations which have the following form~\cite{jackson2008modern}:
\begin{align}
  \label{eq:maxwell_strong}
  \nabla \times \left( \mu_r^{-1} \nabla \times \vec{E} \right) - k_0^2 \varepsilon_r \vec{E} = 0.
\end{align}
Here, $\vec{E}(\vec{x})$ is the complex electric field vector, $\mu_r$ is the relative magnetic permeability, $k_0$ is the free space wave number, and $\varepsilon_r$ is the relative permettivity (may be complex).

We are interested in simulating radio frequency devices that have metallic walls, lossy dielectric materials, open boundaries, and excitations in the form of lumped ports or waveguide ports. The following boundary conditions (BC) are used:
\begin{align}
  &\text{perfect electric conductor (PEC)}
  &&\vec{n} \times \vec{E} = 0,
  \\%
  &\text{first order absorbing BC (ABC)}%
  &&\vec{n} \times \nabla \times \vec{E} = - j k_0 \vec{n} \times \vec{n} \times \vec{E}_t,
  \label{eq:ABC}\\%
  &\text{rectangular lumped port BC}%
  &&\vec{n} \times \nabla \times \vec{E} 
  = -j k_0 \frac{Z_0}{Z_s} \vec{n} \times \vec{n} \times (\vec{E} - 2 \vec{E}_0),
  \label{eq:lumped_port_bc}%
\end{align} 
where $\vec{n}$ is the outward unit normal vector for the boundary, $Z_0$ is the free space impedance, and $\vec{E}_0$ is the excitation electric field.

In \cref{eq:lumped_port_bc}, rectangular lumped ports are bounded by two parallel PEC edges and two parallel perfect magnetic conductor (PMC) edges. For such ports, the electric field $\vec{E}_0$ is uniform and directed from PEC to PEC, and $Z_s = ZW/L$, where $Z$ is the port resistance in Ohm, $W$ is the length of the PEC edge, and $L$ is the length of the PMC edge.

For brevity, we provide expressions only for the lumped port, the coaxial port with one transverse electromagnetic mode is treated similarly.

\section{Finite Element Method and a~Posteriori Error Estimates}%
\label{sec:fem}

We employ the standard finite element method with second order Nédélec edge basis functions of the first kind~\cite{Nedelec}. The weak formulation of the problem has the following form (for lumped ports):
\begin{multline*}
    \int_{\Omega}\left( \nabla \times \vec{\phi} \cdot \left(\mu_r^{-1}\nabla \times \vec{E}\right)- k_0^2 \varepsilon_r \vec{\phi} \cdot  \vec{E}\right) \,\D V
    +\int_{S_{ABC}} \frac{j k_0}{\mu_r} \left( \vec{n} \times \vec{E}\right) \cdot \left( \vec{n} \times \vec{\phi}\right)  \,\D S
    \\
    + \sum_i \int_{S_{P_i}} \frac{j k_0 Z_0}{\mu_r Z_{s,i}}\left( \vec{n} \times \vec{E}\right) \cdot \left( \vec{n} \times \vec{\phi}\right)  \,\D S
    = \int_{S_{P_m}} \frac{2j k_0 Z_0}{\mu_r Z_{s,m}} \vec{n} \times \vec{E}_0 \cdot \vec{n} \times \vec{\phi}  \,\D S
    ,
\end{multline*}
where $S_{ABC}$ is the truncation boundary of the computational domain, $S_{P_i}$ is the surface of port $i$, $Z_{s,i}$ is the sheet impedance of port $i$, $m$ is the index of the excited port, $\vec{E}_{0}$ is the excitation electric field, and $\vec{\phi}$ is a test function.
Here the solution $\vec{E}$ and the test function $\vec{\phi}$ are from the $H(\curl, \Omega)$ functional space, i.e., both $\vec{E}$ and $\curl \vec{E}$ belong to the standard Lebesque space $L_2(\Omega)$, see~\cite{Nedelec} for rigorous definitions.

\runinhead{Error Estimates}
To drive AMR, one needs element-wise error estimates or error indicators. There are many different types of error estimates: explicit residual estimates~\cite{Botha20053717, Monk1998173}, implicit estimates based on solutions of local element-wise or patch-wise problems~\cite{Schoberl2008633}, and adjoint-based goal-oriented error estimates~\cite{Ingelstrom2006, Sun_Cendes_Lee2000}, to name a few.

In our workflow, we use a goal-oriented error indicators similar to~\cite{Sun_Cendes_Lee2000}. The indicators are based on the right-hand side (RHS) of the equations, which is a weighted sum of S-parameters or a multiple of an S-parameter for one excitation. Our experiments have confirmed that these error indicators provide a reasonable trade-off between the accuracy of the estimates and the computational effort to compute them. The error indicator $\eta_T$ for each tetrahedron $T$ is computed using the element and face values $\alpha_T$ and $\alpha_{\Delta, f}$:

\begin{align*}
  \eta_T &= \alpha_T + \sum_{f \in \text{faces$(T)$}} \frac{\alpha_T}{\alpha_T + \alpha_{N_f}} \alpha_{\Delta, f}
  ,\\
  \alpha_T &= - \frac{j \omega \E^{-jk_0 R_T}}{4 \pi \varepsilon R_T} Q_T^2
    - \frac{j \omega \mu \E^{-j k_0 R_T}}{4 \pi R_T} \vec J_T \cdot \vec J_T
  ,\\
  \alpha_{\Delta, f} &= - \frac{j \omega \E^{-jk_0 R_f}}{4 \pi \varepsilon R_f} Q_f^2 
    - \frac{j \omega \mu \E^{-j k_0 R_f}}{4 \pi R_f} \vec J_f \cdot \vec J_f
  ,
\end{align*}
where $N_f$ is the face neighbor of $T$ w.r.t.\ face $f$.
These values are computed from the numerical FEM solution $\vec{E}_h$ using four types of residuals:
\begin{align*}
  Q_T &= \int_T \nabla \cdot (\varepsilon \vec{E}_h) \,\D \vec{x}
  ,\\
  \vec{J}_T &= -\int_T  \frac{1}{j \omega \mu_0} (\nabla \times \mu_r^{-1} \nabla \times \vec{E}_h - k_0^2 \varepsilon_r \vec{E}_h ) \,\D \vec{x}
  ,\\
  Q_f &= \int_f {[\varepsilon \vec{E}_h]}_f \cdot \vec{n} \,\D S
  ,\\
  \vec{J}_f &= -\frac{1}{j \omega \mu_0} \int_f {[\mu_r^{-1} \nabla \times \vec{E}_h]}_f \times \vec{n} \,\D S
  .
\end{align*}
Here ${[\cdot]}_f$ denotes a face jump of some quantity, ${[u]}_f = (u_1 - u_2)$ with normal $\vec{n}$ to the face $f$ directed from side $1$ to side $2$, and $R_T$ and $R_f$ are the tetrahedron and face diameters, respectively.

\Cref{tab:error_indicators} shows how two types of estimators predict tetrahedrons with largest errors. We can order all elements w.r.t.\ to exact errors and error indicators and select sets of $N$ elements with largest errors $top^{\text{exact}}(N)$ and $top^{\text{estim}}(N)$.
To get an approximation to the true error, we computed the solution for the antenna array unit on a very fine uniformly refined mesh with about one million elements and the test mesh and took the $L_2$ norm of the difference between the two solutions.
The table contains values of the ratio
\[
   r(N) = \frac{\abs*{top^{\text{exact}}(N) \cap top^{\text{estim}}(N)}}{N} \cdot \qty{100}{\percent}
   .
\]
which shows how many elements with $N$ largest indicators are among the elements with $N$ largest true errors.

\begin{table}[t]%
  \caption{Fractions of correctly predicted high-error tetrahedra}%
  \label{tab:error_indicators}%
  \centering{}%
  \begin{tabular}{lrrr}%
    \toprule%
    fraction of top-error tetrahedra
      & \qty{ 1}{\percent} & \qty{ 5}{\percent} & \qty{10}{\percent}\\%
    \midrule%
    Indicators~\cite{Sun_Cendes_Lee2000}
      & \qty{26}{\percent} & \qty{44}{\percent} & \qty{54}{\percent}\\%
    Indicators~\cite{Monk1998173}
      & \qty{31}{\percent} & \qty{48}{\percent} & \qty{56}{\percent}\\%
    \bottomrule%
  \end{tabular}%
\end{table}

\begin{figure}[t]%
  \centering{}%
  \includegraphics[width=0.80\textwidth]{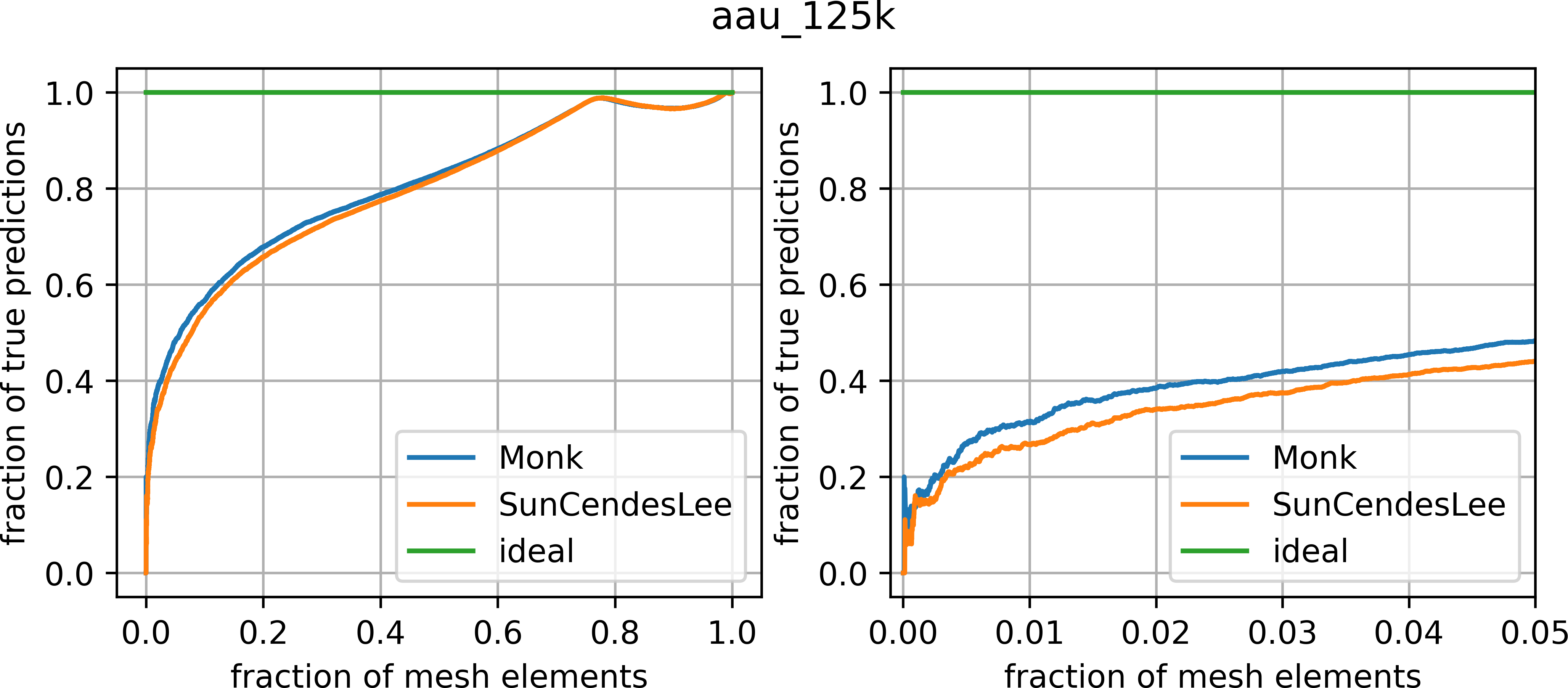}%
  \caption{Accuracy of error indicators for the AAU with \qty{125}{k} mesh elements.
    Right: close-in for the fraction from $0$ to $0.05$}%
  \label{fig:error_indicators}%
\end{figure}

The same ratio is shown in \cref{fig:error_indicators}. It is obvious that both presented indicators are far from perfect, however, they provide a satisfactory quality to drive AMR. Also note that the indicators from~\cite{Sun_Cendes_Lee2000} are slightly worse than the indicators from~\cite{Monk1998173} w.r.t.\ $L_2$ error, which is quite natural since the former indicators are designed for a functional that strongly differs from the $L_2$ error functional.

\section{Tetrahedral Meshing for Models with Thin Material Layers}%
\label{sec:thin_layers_meshing}

The objective of the mesh generation research is the construction of initial tetrahedral meshes and the development of mesh adaptation technology with the aim of achieving adaptive remeshing capabilities similar to state-of-the-art meshing tools, such as those incorporated in the Ansys HFSS industrial software.
The technological chain of the developed meshing tool is based on the following sequence:
\begin{itemize}
\item Analysis of the input 3D model using Open Cascade based representation structures: vertices, edges, volumes.
\item Creation of the initial tessellation using the Open Cascade library.
\item Construction of the initial mesh using an in-house tetrahedral meshing algorithm  similar to the incremental constrained Delaunay triangulation implemented in \tetgen{}~\cite{tetgen}.
\item Mesh optimization using an in-house optimization tool.
\item Mesh adaptation using heuristic finite element a~posteriori error estimates combined with mesh refinement and smoothing.
\end{itemize}

\subsection{Structure of the Underlying Algorithms}%
\label{sec:panel_tet_structure}

The control of the initial tetrahedral mesh generation is based on the following control parameters:
\begin{itemize}
\item \emph{chordal error}: measures distance from tessellation to the surface of the model;
\item \emph{angular resolution}: the deviation of the boundary mesh face normals from the exact normals;
\item \emph{mesh size upper bound}: a constant defined for each material that limits the maximum mesh edge length;
\item \emph{mesh quality bound}: a measure of the admissible mesh cell shape distortion;
\item \emph{mesh anisotropy bound}: a measure of admissible anisotropy applied to the refinement of boundary edges and faces.
\end{itemize}
\emph{Chordal error} and \emph{angular resolution} are responsible for controlling the deviation of the discretized model from the exact one, while \emph{mesh size bound} and \emph{mesh quality bound} are used to control the mesh quality.

We use the following algorithm:
\begin{itemize}
\item The input surface triangulation is used to build an initial tetrahedral mesh with a minimal number of Steiner vertices using the standard Constrained Delaunay Tetrahedralization algorithm~\cite{tetgen}.
\item Then, tetrahedra violating the size and quality constraints are marked for refinement.
We use an algorithm similar to the well-known Ruppert refinement algorithm~\cite{Ruppert} and its generalizations by Shewchuk (see details in~\cite{tetgen}) with some important modifications that prevent over-refinement near thin material layers (see \cref{sec:insertion_strategy}).
\item
   Finally, the mesh from the previous stage is post-processed using both topology optimization (edge/face flips, insertion of virtual boundary layers) and optimization by mesh vertex movement (the smoothing algorithm is described in \cref{sec:smoothing}).
\end{itemize}
It should be noted that smoothing is a crucial part of the algorithm because both initial meshes and meshes resulting from adaptive mesh refinement have slivers with dihedral angles close to $\pi$. These slivers cannot be removed by standard flips and, if not removed, will cause the numerical simulation to fail.

\subsection{New Anisotropic Encroachment Domain (AED) Mesh Vertex Insertion Strategy to Account for Thin Material Layers}%
\label{sec:insertion_strategy}
\begin{figure}[t]%
  \includegraphics[width=1.0\textwidth]{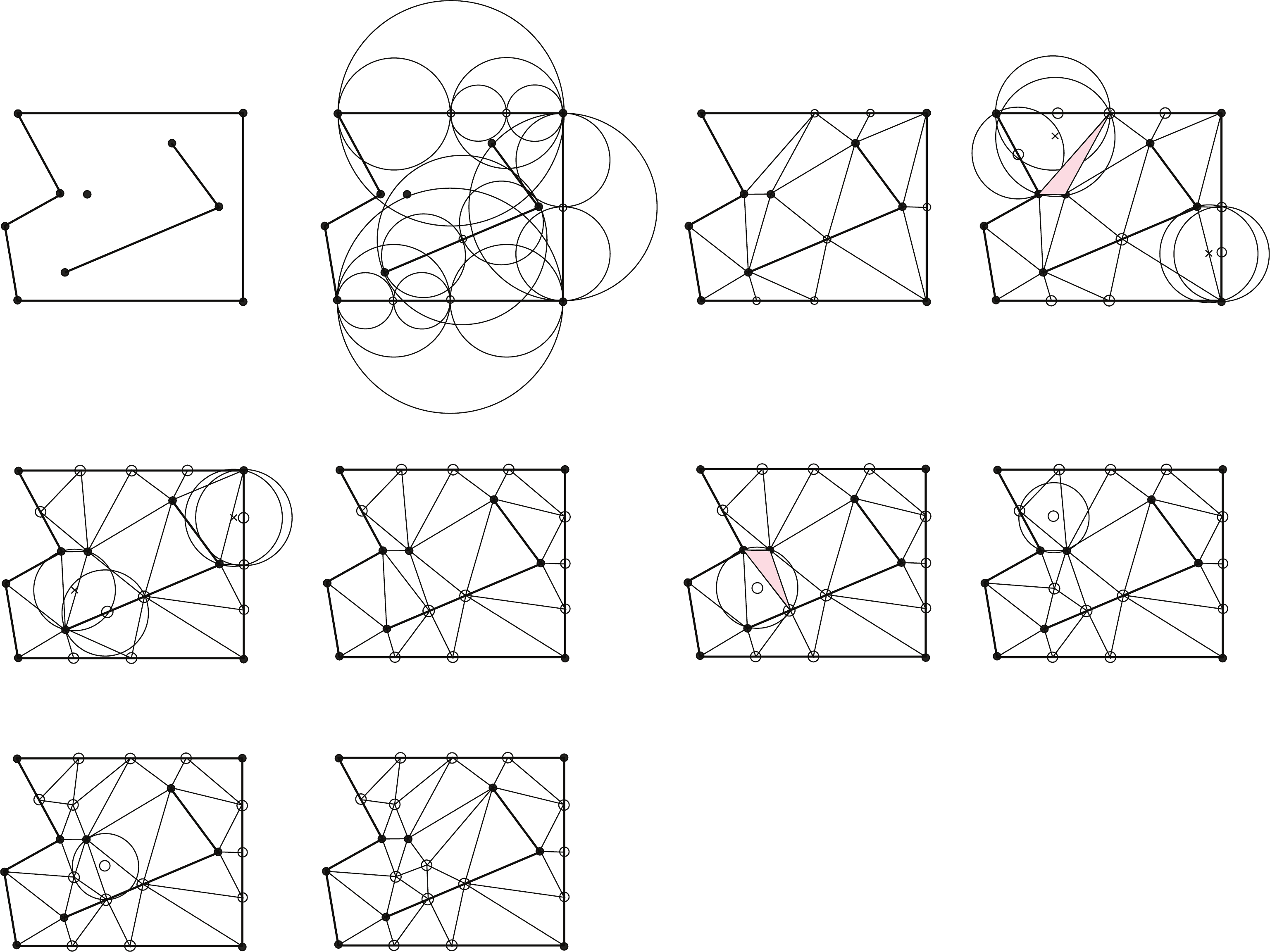}%
  \caption{Ruppert's algorithm simultaneously refines PPLG edges and triangles to meet quality and mesh size requirements while preserving the Delaunay property of the triangulation}%
  \label{fig:Delaunay_refinement_Ruppert}%
\end{figure}
The most appropriate approach for both the initial mesh construction for the FEM solver and for the adaptive refinement procedure is the Delaunay refinement (DR) algorithm. The basic DR algorithm is based on the insertion of the circumcenters of selected mesh cells and specially devised refinement rules for the boundary edges and faces.
\Cref{fig:Delaunay_refinement_Ruppert} provides a simple 2D illustration of the classical Ruppert's DR algorithm and shows the initial splitting of the edges of a planar piecewise linear graph (PPLG) to make their diametral circles empty.
Then, the circumcenters of mesh elements that violate the quality threshold or have a too large circumradius are added as new candidate vertices for the insertion list. If the candidate vertex is located inside the circumcircle of the boundary edge, this vertex is eliminated from the candidate list, while the edge center is added to the list instead. This algorithm is well studied and has quality guarantees under certain assumptions on the initial set of boundary edges~\cite{Ruppert}.

The standard quality measure $Q(T)$ for a triangle $T$ is the ratio of the circumradius to the minimum edge length. It is possible to specify a target quality threshold $Q(T) < Q^*$. Note, however, that even specifying very large values of $Q^*$ does not prevent the generation of very fine meshes in cases imitating thin material layers, e.g., if very closely located long parallel edges are present in the input PPLG.

Similar behavior is observed in 3D. In this case, the input is a triangulated polyhedral complex (PLC). Since surface triangles can be quite bad, it is generally not always possible to obtain Delaunay or constrained Delaunay tetrahedral meshes due to complicated boundary configurations. Therefore, a local violation of the Delaunay property is allowed.

In 3D, similar point generation rules apply: diametrical spheres of both segments and boundary triangles must be empty. To achieve this property, many boundary segments and faces must be split, even if all tetrahedra in the mesh satisfy the size and quality requirements. This can lead to a very strong over-refinement for geometries with high aspect ratio elements.
Experiments with real-life structures such as cavity resonators and antenna elements have confirmed that the standard boundary treatment~\cite{tetgen} leads to over-refinement of the domain boundaries and produces meshes that are too dense.

One possible way to overcome this problem is to build over-refined meshes using the basic DR algorithm and then apply mesh decimation algorithms to attain the prescribed mesh size while respecting the given mesh quality and anisotropy threshold. This algorithm is quite robust but can be quite slow for large and complex models.

To solve the same problem without the over-refinement decimation stage, we propose a modification of the mesh refinement rules near boundaries. The resulting algorithm is:
\begin{itemize}
\item for the tetrahedral mesh $\mathcal{T}$, create the set $T_r$ of candidate tetrahedra marked for refinement;
\item for tetrahedra from the set $T_r$, create ``centers'' which are added to the list of candidate vertices for insertion;
\item create anisotropic encroachment domains (AED) around boundary edges and faces;
\item if the candidate vertex is inside the AED, it is excluded from the list of candidates, while the ``center'' of the encroachment domain is computed and added to the list of candidate vertices instead;
\item add new vertices to the existing mesh;
\item use local mesh optimization by edge flipping and local vertex relocation.
\end{itemize}

\begin{figure}[t]%
  \centering{}%
  \includegraphics[width=1.0\textwidth]{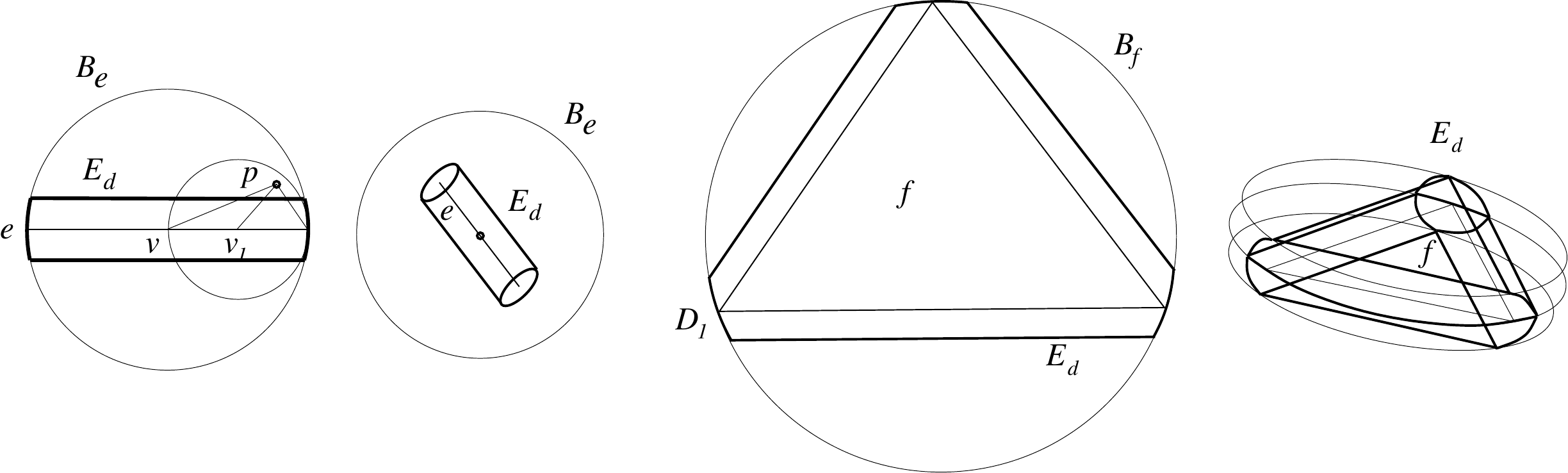}%
  \caption{Anisotropic encroachment domains for boundary edges and faces}%
  \label{fig:Delaunay_refinement_ED}%
\end{figure}

Let us explain the concept of the anisotropic encroachment domain using \cref{fig:Delaunay_refinement_ED} as an illustration. Suppose that the anisotropy parameter $A > 1$ is prescribed. For an arbitrary boundary edge $e$ we assign a diametral circumball $B_e$ with radius $R_e$ and consider the set $C_e$ of the 3D points that are at a distance not exceeding $R_e / A$ from $e$. The resulting figure is obviously a rounded cylinder. We define the encroachment domain $E_d(e)$ as the intersection of the rounded cylinder $C_e$ and the circumball $B_e$,
\[
E_d(e) = C_e \cap B_e.
\]
For the boundary face $f$, we first construct an equatorial circumball $B_f$ with radius equal to $R_f$ and define the closest set $C_f$ of the 3D points that are at a distance not exceeding $R_e / A$ from $f$. The resulting figure resembles a rounded triangular prism with the thickness equal to $2 R_e / A$. The face encroachment domain $E_d(e)$ is defined as the intersection of the rounded prism $C_f$ and the face circumball $B_f$,
\[
E_d(e) = C_f \cap B_f.
\]

In \cref{fig:Delaunay_refinement_ED} (left), trying to insert the candidate vertex $p$ using the standard Ruppert refinement based on the diametral circumballs around boundary edges would force splitting the edge $e$ into three sub-edges. In the proposed scheme with $A=5$, vertex $p$ can be inserted into the mesh without splitting the edge.

If the insertion center $c$ falls inside $E_d(f)$ of a certain boundary face $f$, one should first check if it is inside the AED of the edges of $f$. If so, the edge processing branch is used. Note that this simple modification essentially suppresses attempts to create a new boundary vertex outside of its generated face, while allowing a controlled level of distortion.

\subsection{Pre-processing and Post-Optimization for Sliver Elimination}

The appearance of special flat or nearly flat tetrahedra called ``slivers'' is one of the most annoying features of the Delaunay tetrahedral meshers.
Let us explain this important concept in more detail.

\begin{figure}[t]%
  \centering{}%
  \subcaptionbox{\label{fig:sliver-definition:a}}{%
    \includegraphics[height=0.15\textheight]{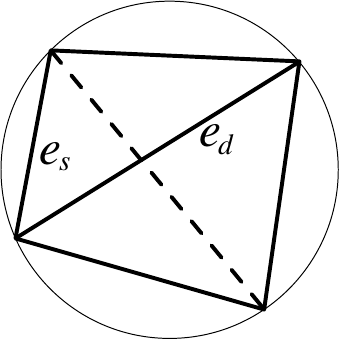}%
  }%
  \hfill{}%
  \subcaptionbox{\label{fig:sliver-definition:c}}{%
    \includegraphics[height=0.15\textheight]{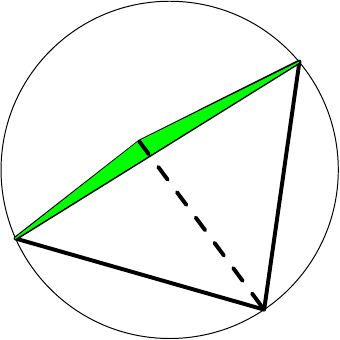}%
  }%
  \hfill{}%
  \subcaptionbox{\label{fig:sliver-definition:b}}{%
    \includegraphics[height=0.15\textheight]{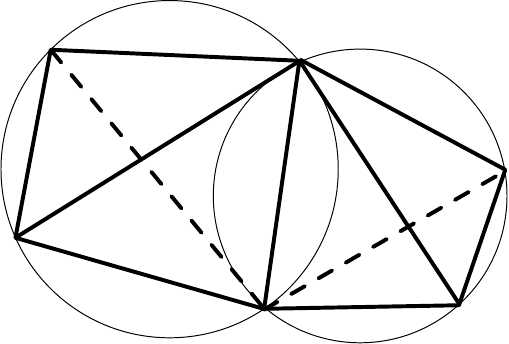}%
  }%
  \\%
  \subcaptionbox{\label{fig:sliver-definition:d}}{%
    \includegraphics[height=0.15\textheight]{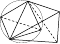}%
  }%
  \hspace{0.12\textwidth}%
  \subcaptionbox{\label{fig:sliver-definition:e}}{%
    \includegraphics[height=0.15\textheight]{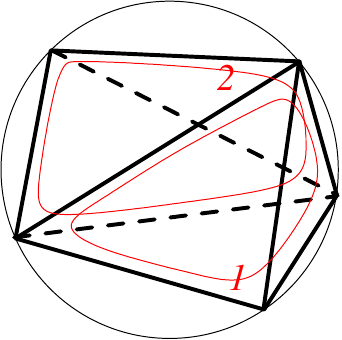}%
  }%
  \caption{Slivers:
    \protect\subref{fig:sliver-definition:a} sliver with side and diagonal edges,
    \protect\subref{fig:sliver-definition:c} constrained (blocked) sliver,
    \protect\subref{fig:sliver-definition:b} weak coupling of slivers via side edges,
    \protect\subref{fig:sliver-definition:d} strong coupling of slivers via diagonal edges, and
    \protect\subref{fig:sliver-definition:e} triangulation of a Delaunay face into nearly flat tetrahedra results in sliver chains}%
  \label{fig:sliver-definition}%
\end{figure}

As shown in \cref{fig:sliver-definition:a}, a sliver is a flat or nearly flat tetrahedron with four ``side'' edges $e_s$ with nearly zero dihedral angles and two ``diagonal'' edges $e_d$ with dihedral angles close to $\pi$.
Slivers are very difficult to visualize and eliminate: the circumradius of a sliver is not large, so that the standard quality criteria based on the ratio of the circumradius to the minimum edge length cannot be used to identify and eliminate them.
The main sources of slivers are co-circular point sets, which would normally produce flat polygonal Delaunay faces, but due to the rounding errors may or may not produce almost flat tetrahedra instead.
Most of these Delaunay faces are quadrilaterals that produce isolated slivers.

Boundary effects lead to more complicated cases. By definition, a constrained Delaunay tetrahedron has an empty \emph{visible} part inside the circumball, hence almost flat tetrahedra can appear, similar to the one shown in \cref{fig:sliver-definition:c}. The circumball of such a tetrahedron can be extremely large, but only a small part of this ball is visible from inside the tetrahedron itself. \Cref{fig:sliver-definition:c} shows an example of a ``blocked sliver'' with all four vertices and one face on the fixed boundary.

The majority of isolated slivers can be eliminated by standard 2-to-3 edge flip operations. Some slivers can be eliminated by insertion of their circumcenters into the mesh. In addition, many special methods for sliver removal have been developed, see~\cite{tournois2009perturbing} for reference.  The situation when slivers are coupled by side edges, as in  \cref{fig:sliver-definition:b}, is not critical as well.

The most difficult case for analysis and correction is related to ``sliver chains'', where a coupled pair of slivers has at least one common diagonal edge, as in \cref{fig:sliver-definition:d}. The origin of a ``sliver chain'' is illustrated in \cref{fig:sliver-definition:e}, where a pentagonal Delaunay face is triangulated as a chain consisting of two slivers.

In a simulation, slivers can cause several problems:
\begin{itemize}
  \item large numerical error in FEM interpolation and solution;
  \item high condition number of the FEM linear system and, as a result, slow convergence of the iterative solvers;
  \item failure of direct solvers due to numerical zero pivots and similar problems.
\end{itemize}

To improve the ability of the optimization algorithm to eliminate slivers, one should take care during the pre-processing stage and eliminate the most vexing input configurations. It is very useful to eliminate very obtuse surface triangles, as well as to identify potential polygonal Delaunay faces in the input data and treat them separately.

\runinhead{Elimination of Very Obtuse Triangles on the Surface}%
\label{sec:obtuse}

Tessellations of curved surfaces tend to produce a number of poorly shaped triangular facets. While needle-type elements are essential because they conform well to surfaces with anisotropic curvature or developable surfaces, the obtuse triangles are the source of large finite element simulation errors and potential generator of slivers in the volume meshes. Therefore, eliminating most of the obtuse surface triangles improves both surface and volume meshes. We have developed an algorithm for eliminating obtuse triangles on manifold surfaces. The non-manifold version is still under development.

The idea of the algorithm is the following: we start from the vertex of a triangle with an obtuse angle and split this triangle by a segment, creating a new vertex on the longest edge opposite to the obtuse angle. This new vertex is either connected to a certain existing vertex or the splitting is continued, creating the shortest polylines on the triangulated surface. We use special criteria to finalize each polyline by pointing to an existing vertex or to a vertex already created by another polyline.

Note that the actual splitting of the surface triangles in this algorithm is not done, only the set of new vertices on the edges is created. Then, for each triangle we build an independent Delaunay triangulation to create new split triangles. This algorithm can be repeated several times until all obtuse triangles above a given threshold are exhausted. We have found that setting the threshold below $150$ degrees does not significantly increase the number of triangles and edges in realistic tessellations.
 
\begin{figure}[p]%
  \centering{}%
  \includegraphics[width=0.68\textwidth]{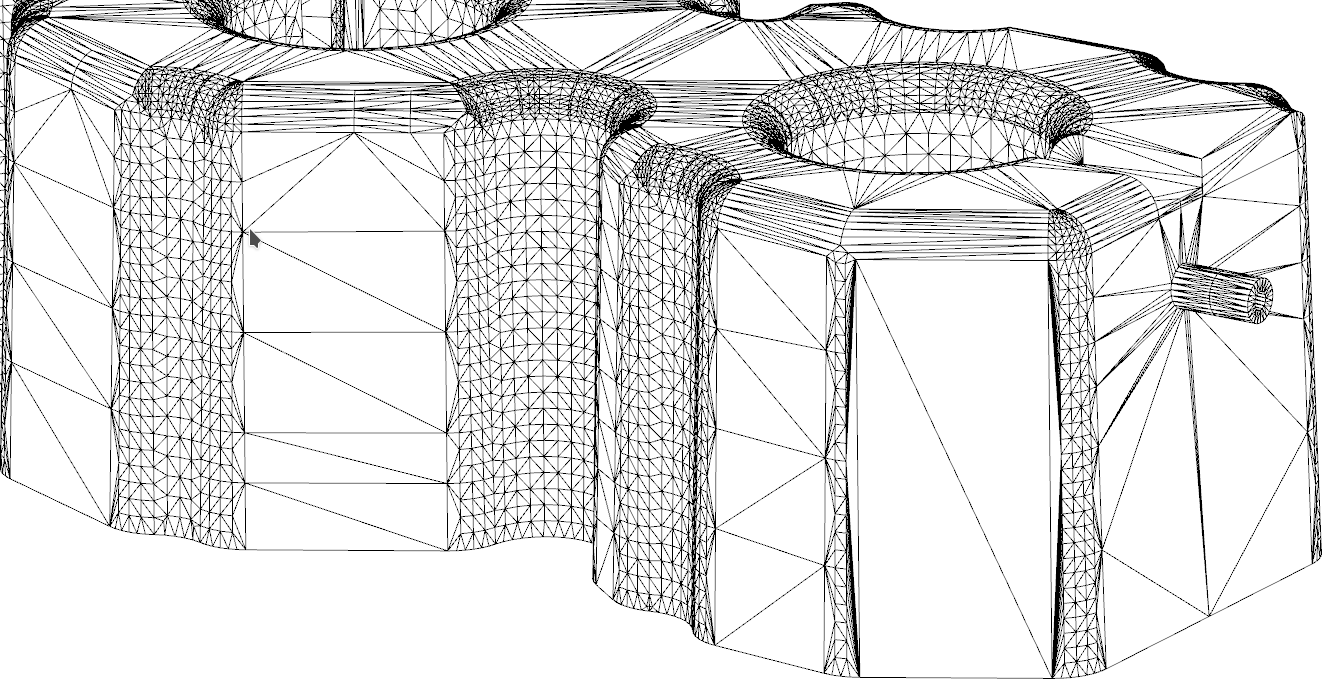}%
  \\[2ex]%
  \includegraphics[width=0.68\textwidth]{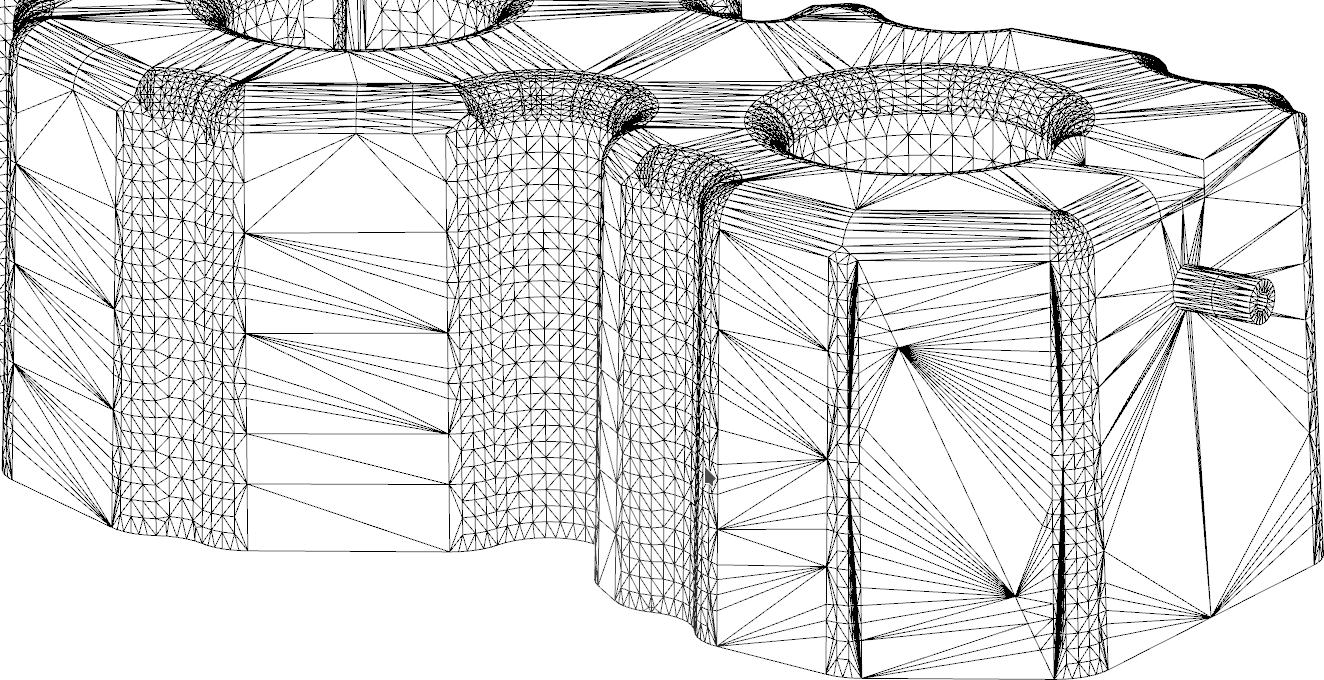}%
 \caption{%
    Top: surface tessellation with obtuse triangles.
    Bottom: Refined tessellation with maximum angle below $150$ degrees%
  }%
  \label{fig:obtuse-1}%
\end{figure}

\begin{figure}[p]%
  \centering{}%
  \includegraphics[height=0.20\textheight]{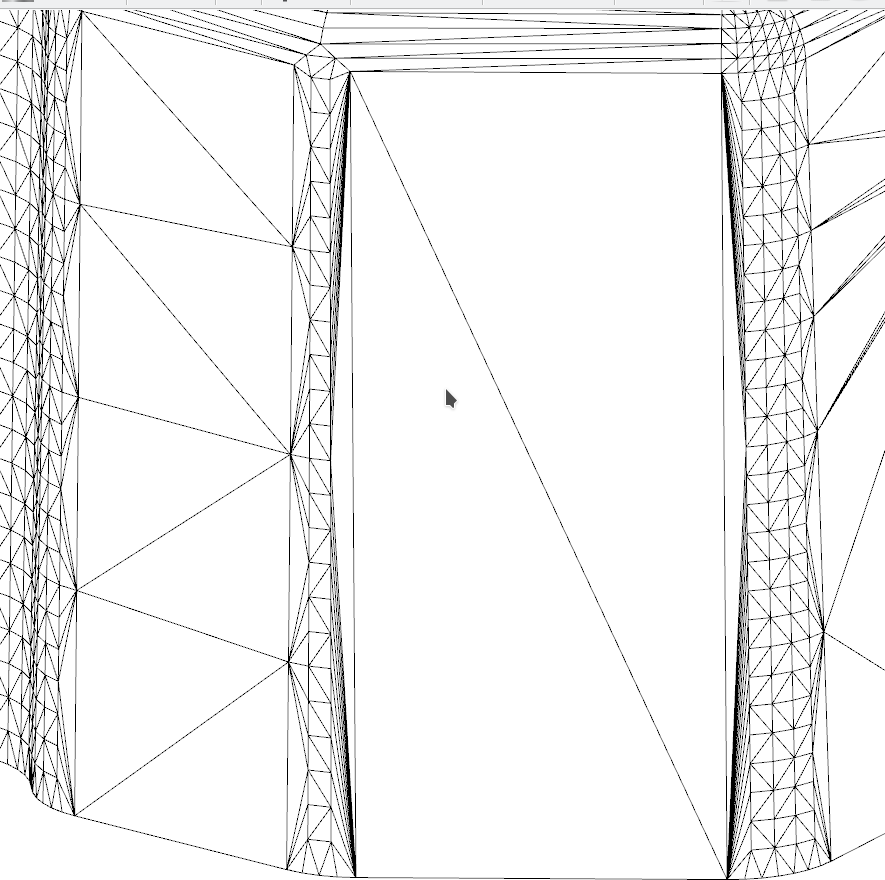}%
  \hspace{0.1\textwidth}%
  \includegraphics[height=0.20\textheight]{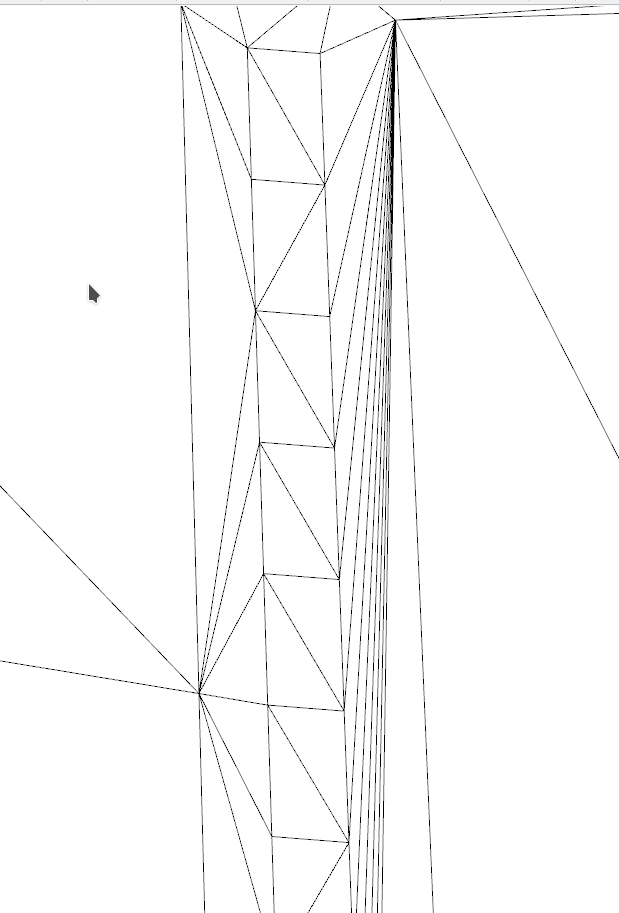}%
  \hspace{0.1\textwidth}%
  \includegraphics[height=0.20\textheight]{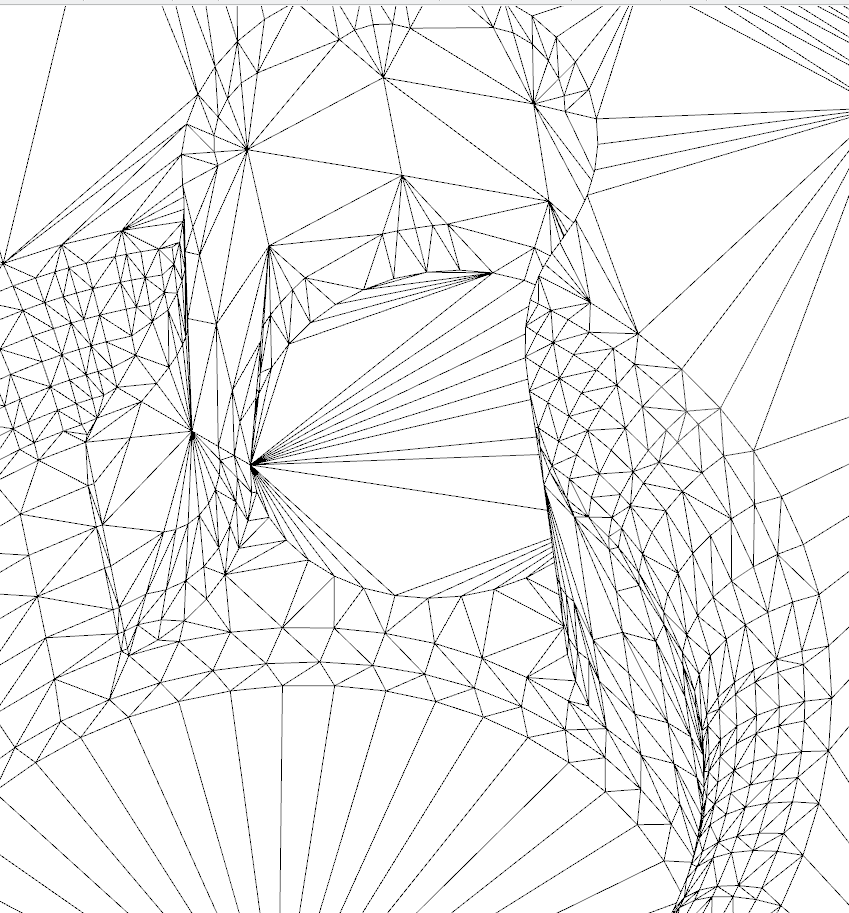}%
  \\[2ex]%
  \includegraphics[height=0.20\textheight]{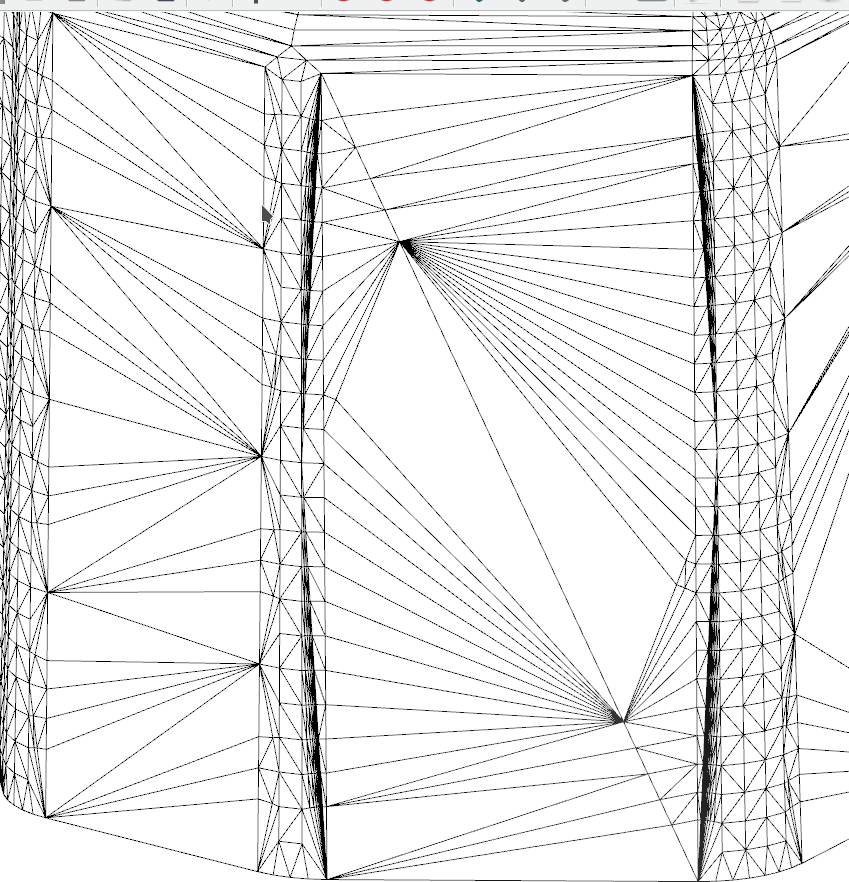}%
  \hspace{0.1\textwidth}%
  \includegraphics[height=0.20\textheight]{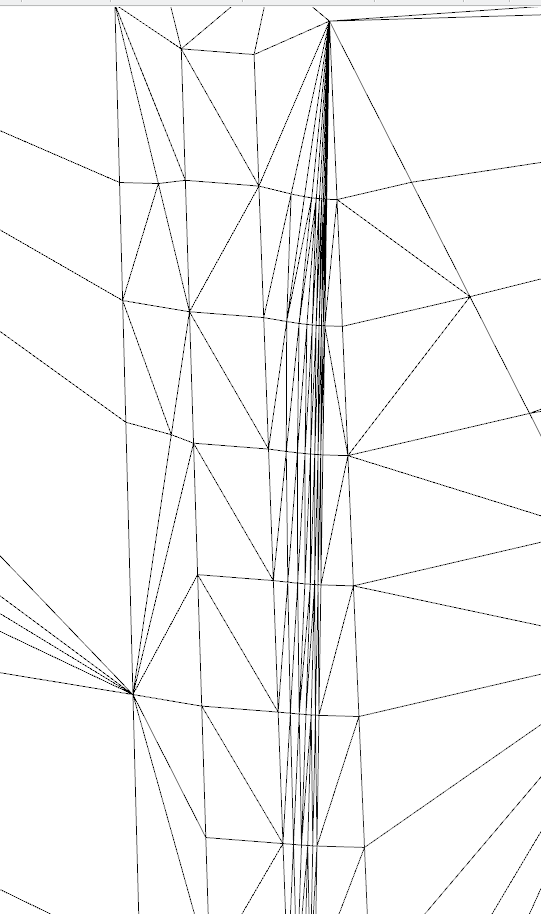}%
  \hspace{0.1\textwidth}%
  \includegraphics[height=0.20\textheight]{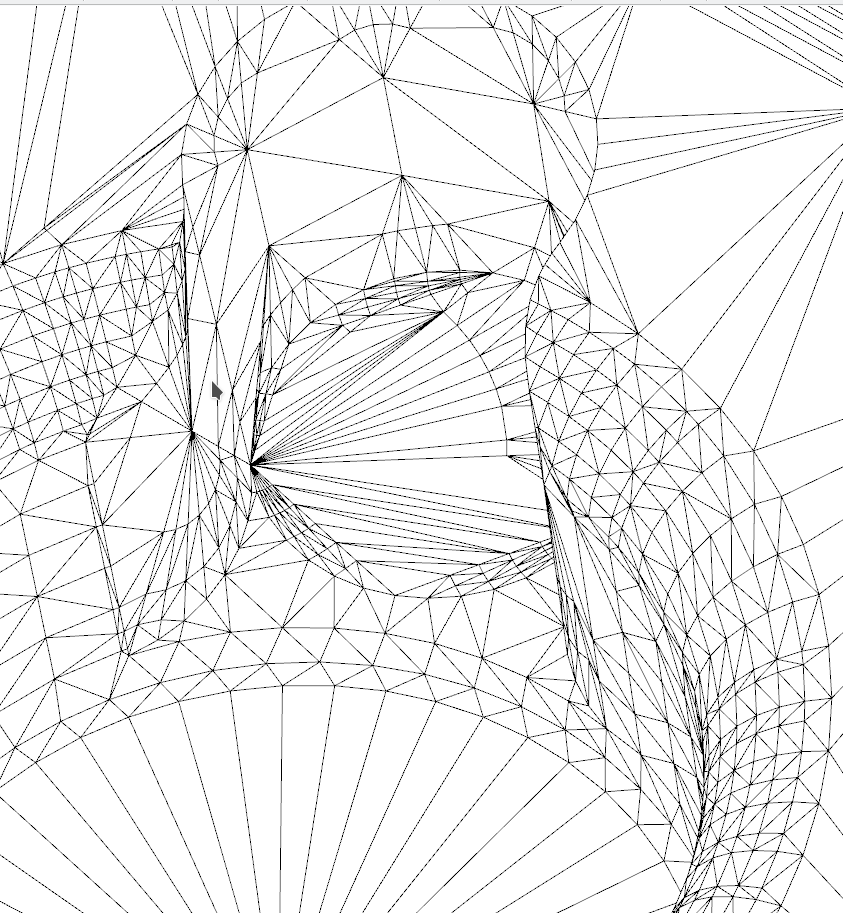}%
  \caption{Enlarged tessellation fragments before (top row) and after (bottom row) refinement}%
  \label{fig:obtuse-2}%
\end{figure}

In the particular case shown in \cref{fig:obtuse-1,fig:obtuse-2} the number of vertices increases from $6491$ to $7516$, while the maximum angle is reduced from $179.9$ to $150$ degrees.

Elimination of obtuse triangles on the surface
helps to make the mesh generation problem less stiff, but it is not enough.

\runinhead{Detection of Polygonal Delaunay Faces}%

\begin{figure}[t]%
  \subcaptionbox{\label{fig:polygonal-slivers:a}}{%
    \includegraphics[height=0.122\textheight]{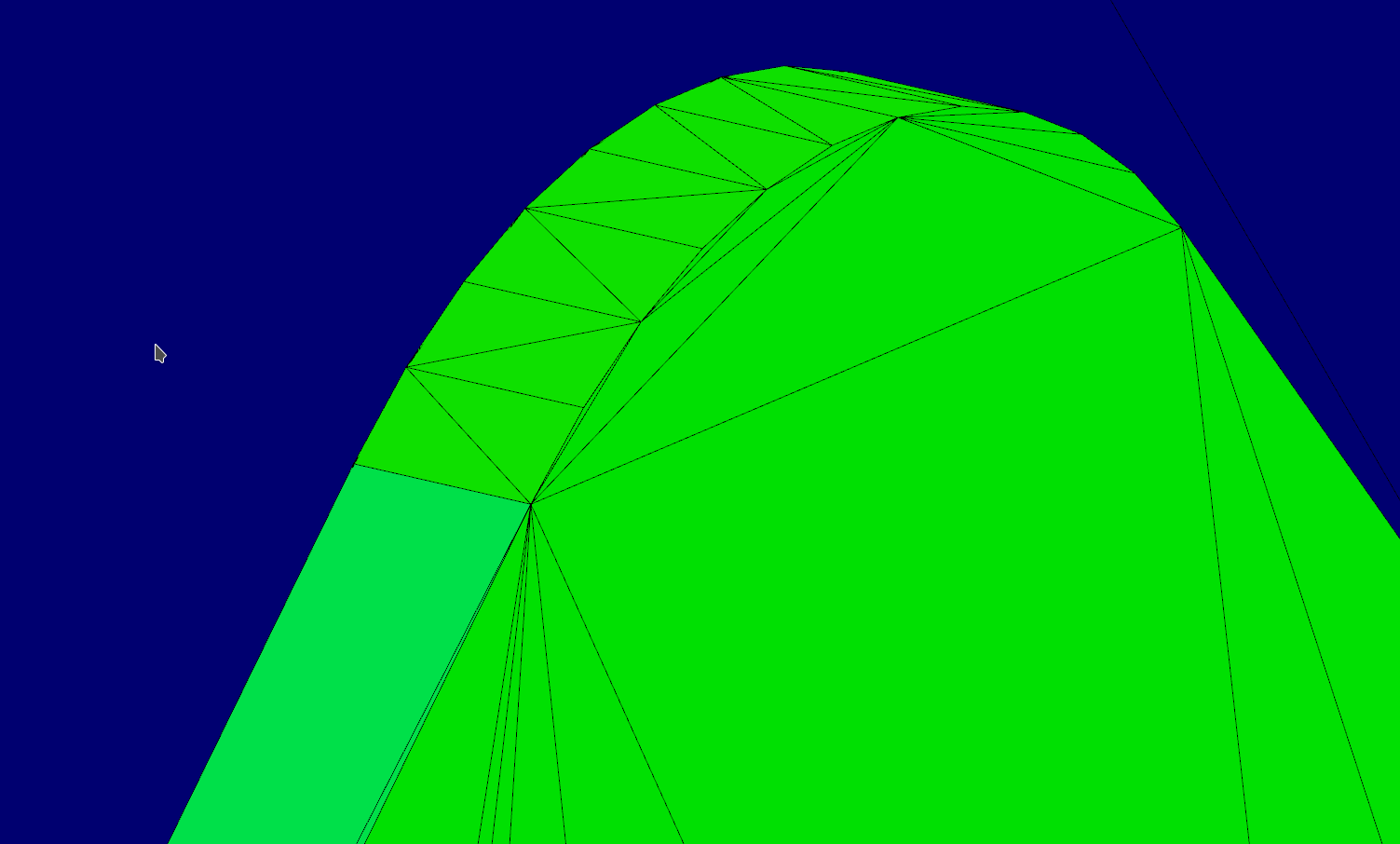}%
  }%
  \hfill{}%
  \subcaptionbox{\label{fig:polygonal-slivers:b}}{%
    \includegraphics[height=0.122\textheight]{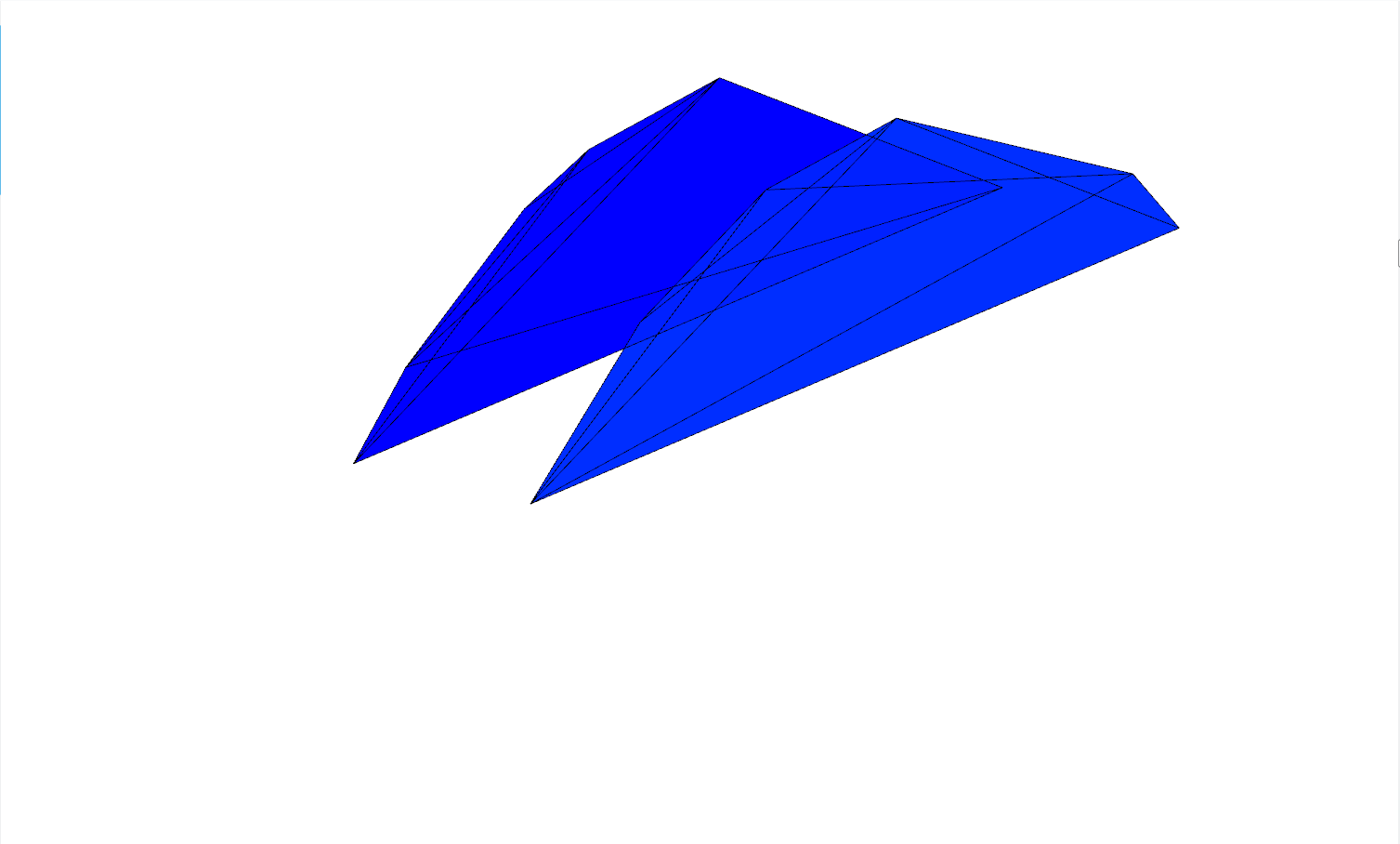}%
  }%
  \hfill{}%
  \subcaptionbox{\label{fig:polygonal-slivers:c}}{%
    \includegraphics[height=0.122\textheight]{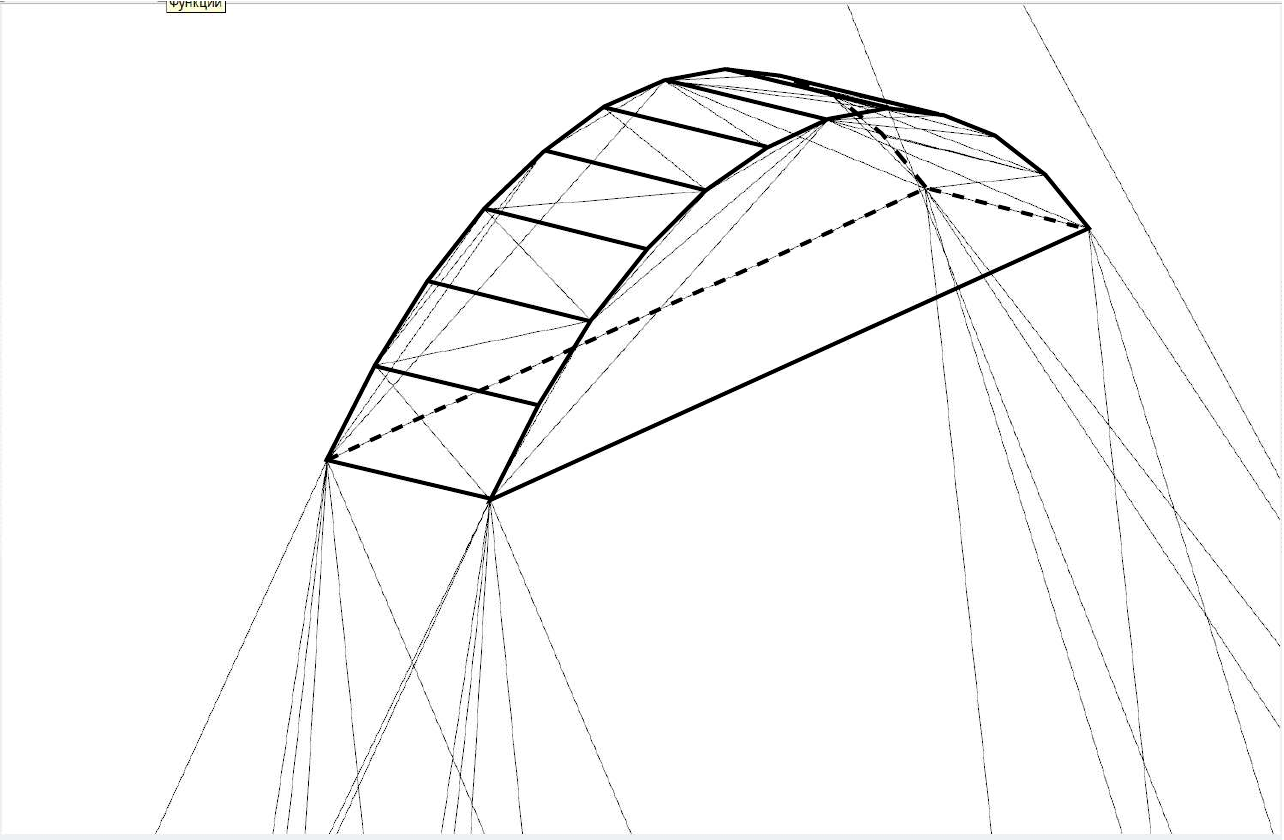}%
  }%
  \caption{Sliver chains:
  \protect\subref{fig:polygonal-slivers:a}~fragment of a tessellation,
  \protect\subref{fig:polygonal-slivers:b}~sliver chains in the mesh,
  \protect\subref{fig:polygonal-slivers:c}~the~corresponding Delaunay polyhedron}%
  \label{fig:polygonal-slivers}%
\end{figure}

A convex envelope of the three-dimensional set of vertices lying on the empty sphere is called a Delaunay polyhedron. The decomposition of a Delaunay polyhedron into tetrahedra is not unique. Note that all tetrahedra in this additional decomposition have the same circumcenter and circumradius. If two adjacent Delaunay polyhedra have inconsistent triangulations, it results in slivers instead of common faces. The following engineering example below shows that this is not an abstract concept. Consider a tessellated fragment of a rib with circular rounding, as in \cref{fig:polygonal-slivers:a}. \Cref{fig:polygonal-slivers:b} shows the resulting sliver chains, and \cref{fig:polygonal-slivers:c} shows the corresponding Delaunay polyhedron which is essentially a polygonal prism.

The proposed solution is to consider the tetrahedral mesh immediately after boundary recovery, which is quite close to the constrained Delaunay mesh, create all Delaunay polyhedra sharing true polygonal faces by constructing clusters of tetrahedra with nearly coinciding circumspheres, and use a consistent decomposition of these polyhedra into a tetrahedral mesh. This is ongoing research.

\runinhead{Variational Mesh Smoothing with Sliding Boundary Vertices}%
\label{sec:smoothing}

The technique for sliver elimination via mesh point movement is based on the concept of a mesh as a deformed elastic body.
Let $\xi_1$, $\xi_2$, $\xi_3$ denote the Lagrangian coordinates associated with the elastic material, and  $x_1$, $x_2$, $x_3$ the Eulerian coordinates of a
material point.  The spatial mapping $x(\xi) : \R^3 \to \R^3$ defines
a stationary elastic deformation, its Jacobian matrix is denoted by $C$, $C_{ij} = {\p x_i}/{\p \xi_j}$.
We search for the elastic deformation  $x(\xi)$ that minimizes the
stored energy functional~\cite{Kudryavtseva-2014}
\begin{equation} \label{weighted-stored-energy-functional}
   F(x) = \int\limits_{\Omega_\xi} {W}(C) \,\D\xi,
\end{equation}
where ${W}(C)$ is the polyconvex elastic potential (internal energy), which is a weighted sum of the shape distortion measure and the volume distortion measure~\cite{Garanzha-2000},
\begin{equation}%
  \label{eq.geometric-elastic-potential}%
  W(C)  =  (1 - \theta) \frac{\frac{1}{3} \tr\left(C^T C\right)}{\det C^{\frac{2}{3}}}
  + \frac{1}{2} \theta \left( \frac1{\det C} +  \det C \right)
.
\end{equation}
In most cases we set $\theta = 4/5$. However, if the size/volume distribution function is not known a~priori, it is more convenient to set $\theta = 0$.

Suppose the domain $\Omega_\xi$ can be partitioned into simplicial cells $U_k$. Then the stored energy functional \cref{weighted-stored-energy-functional} can be approximated by the discrete functional
\begin{equation}%
  \label{half-discret_stored-energy-functional}%
  F (x_h(\xi)) = \sum\limits_{k}  W(C_k) \vol(U_k),
\end{equation}
where  $C_k = \grad x_h(\xi)$ is the Jacobian matrix of the continuous piecewise linear deformation $x_h(\xi)$, which is constant on the $k$-th simplex. To optimize a simplicial mesh with vertex movement along straight boundary edges and flat boundary panels, we use the efficient optimization algorithm developed in~\cite{Kudryavtseva-2014}. We apply a special aggregation procedure, which extracts flat polygonal patches from boundary triangulation and finds sets of straight edges bounding these patches. Essentially, a polyhedral complex defining a discrete model is recovered, which allows us to divide all boundary vertices into fixed vertices, vertices moving along edges, and vertices moving along flat faces.

\begin{figure}[p]%
  \centering{}%
  \subcaptionbox{\label{fig:kubik_smoothing:cube}}{%
    \includegraphics[width = 0.52\textwidth]{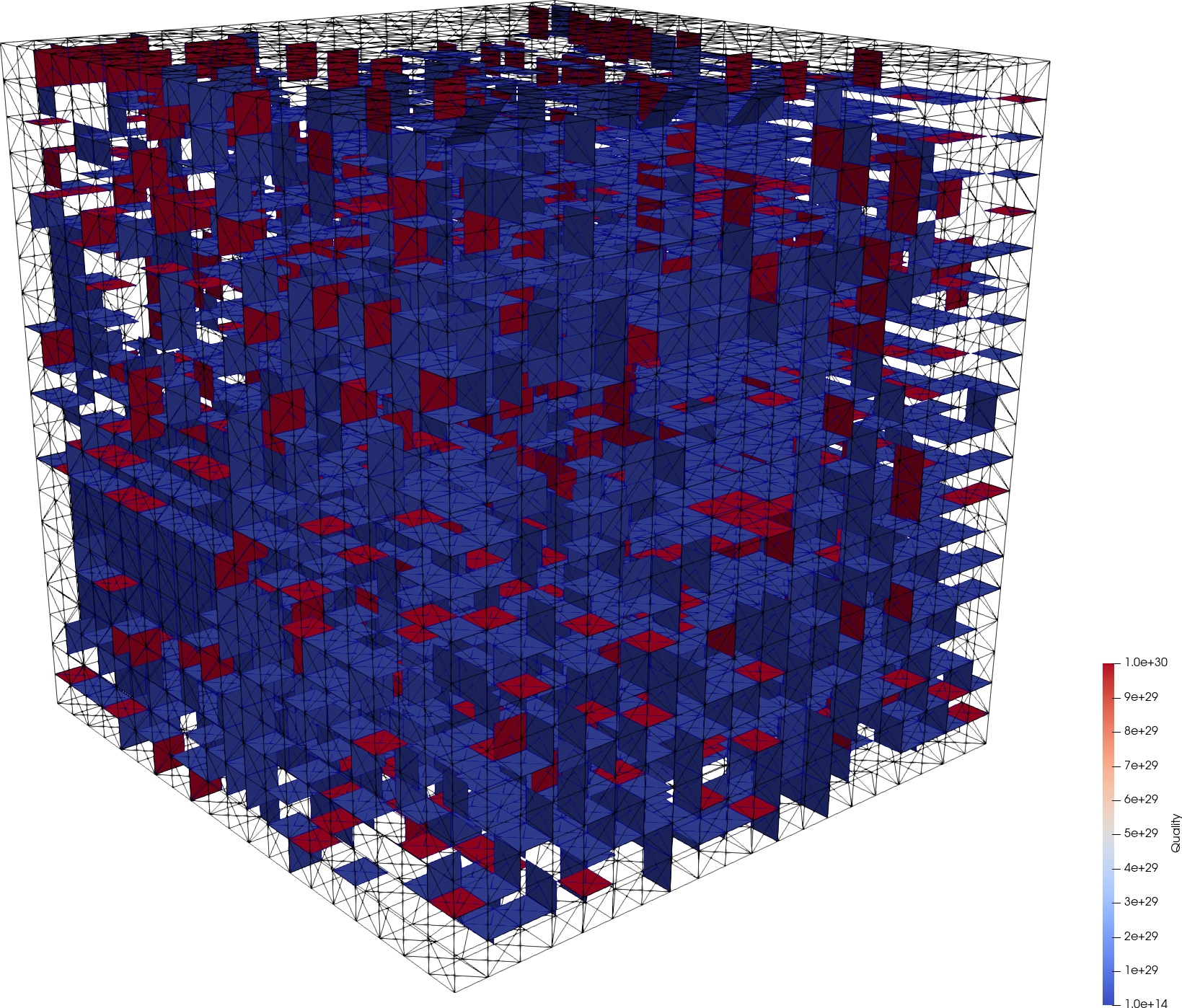}%
  }%
  \hfill{}%
  \subcaptionbox{\label{fig:kubik_smoothing:optimized}}{%
    \includegraphics[width = 0.467\textwidth]{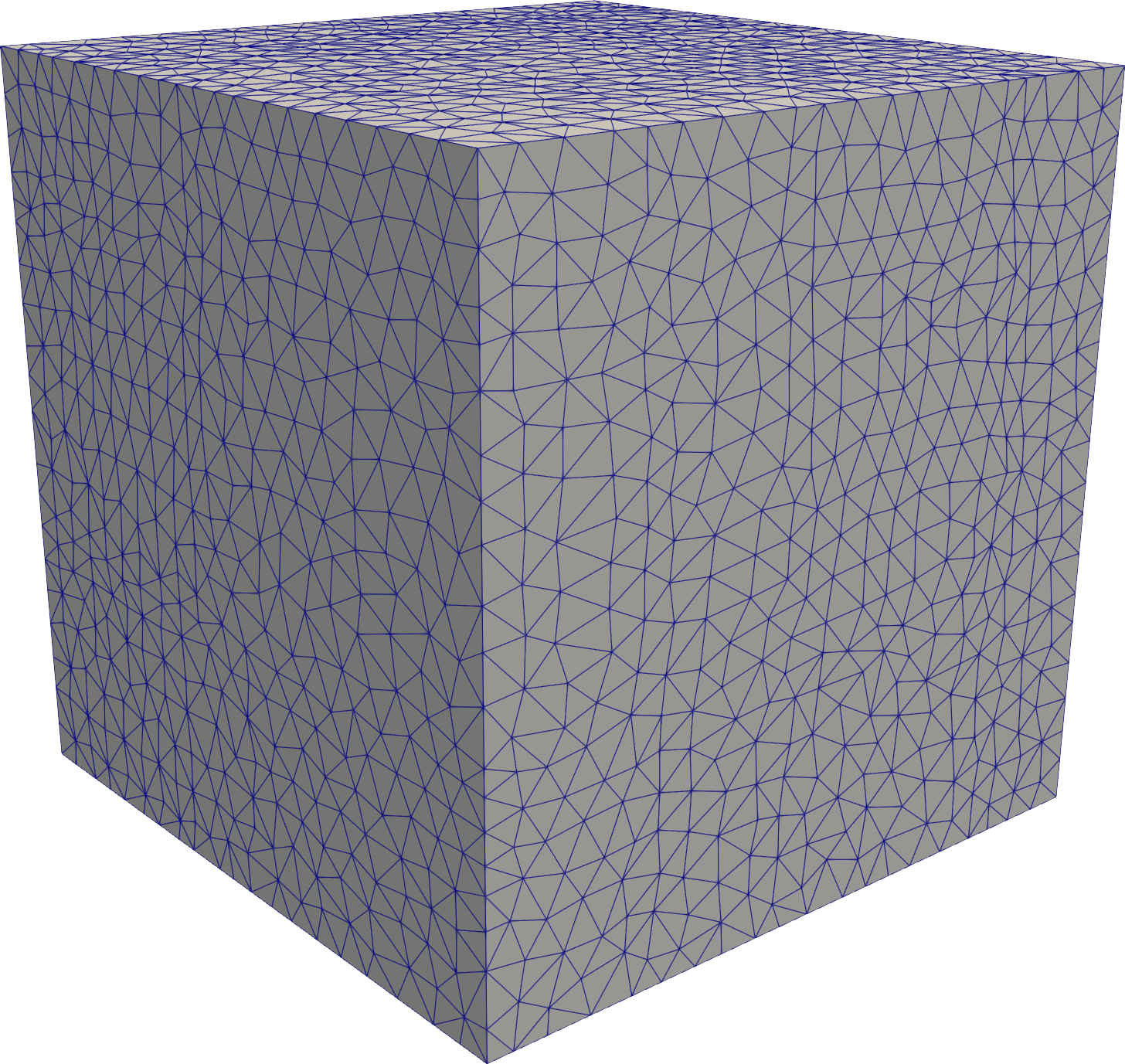}%
  }%
  \\%
  \subcaptionbox{\label{fig:kubik_smoothing:slivers}}{%
  \includegraphics[width = 0.27\textwidth]{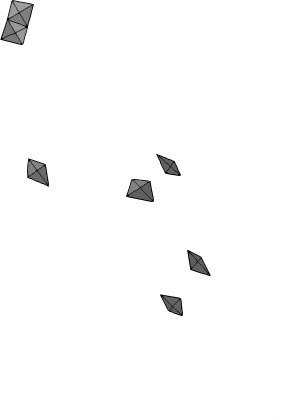}%
  }%
  \hfill{}%
  \subcaptionbox{\label{fig:kubik_smoothing:histo}}{%
    \includegraphics[width = 0.69\textwidth]{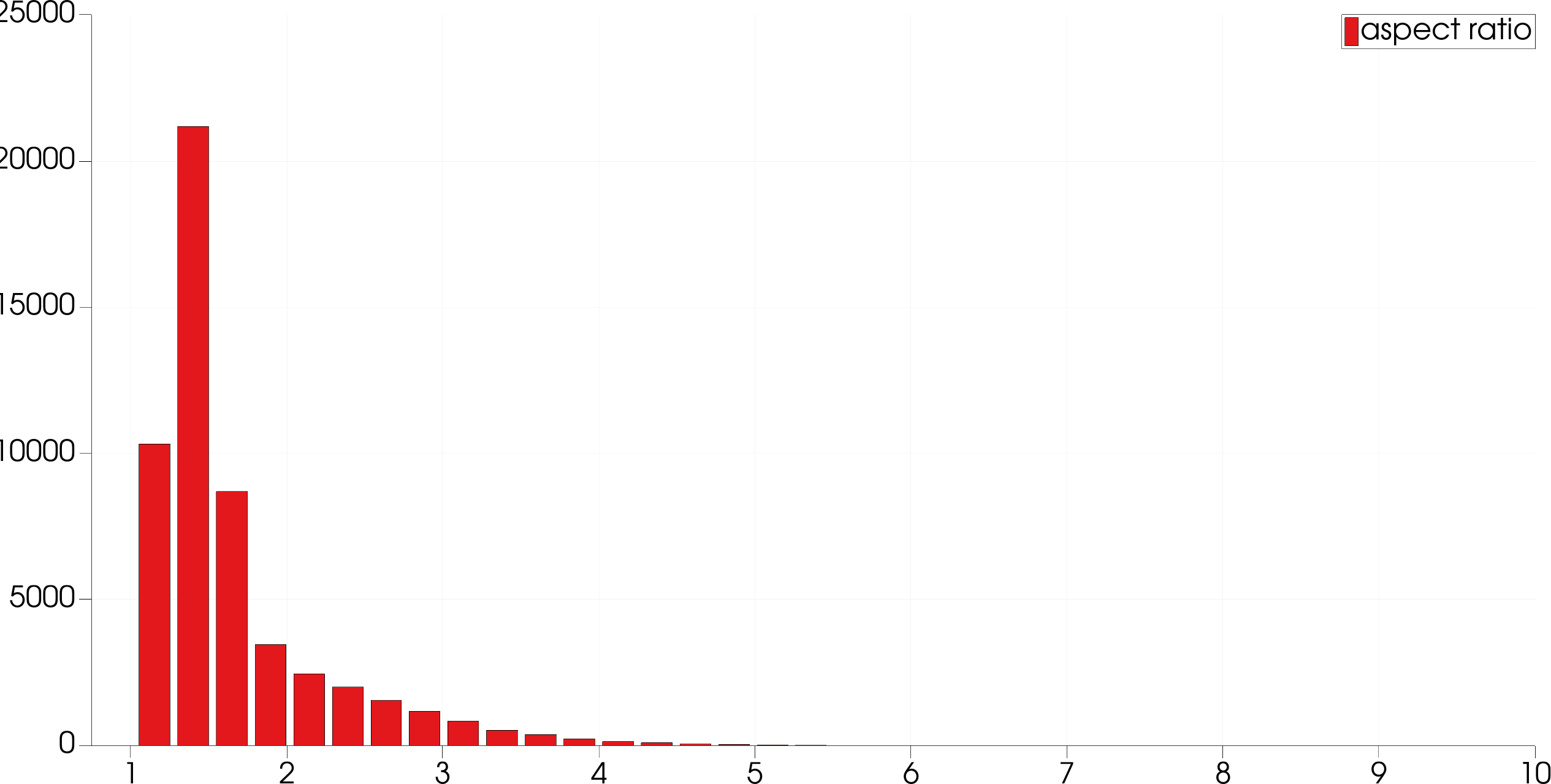}%
  }%
  \caption{Smoothing example:
     \protect\subref{fig:kubik_smoothing:cube} \tetgen{} mesh of a cubic lattice with slivers attached to the faces of the regular cubic mesh,
    \protect\subref{fig:kubik_smoothing:optimized} optimized mesh,
    \protect\subref{fig:kubik_smoothing:slivers} remaining slivers after optimization, and
    \protect\subref{fig:kubik_smoothing:histo} aspect ration histogram of the optimized mesh
  }%
  \label{fig:kubik_smoothing}%
\end{figure}

The following simple test case illustrates the efficiency of the optimization tool. Consider a tetrahedral mesh for a set of vertices of a cubic lattice under arbitrary rotation introducing rounding errors in the input data. For such an input, \tetgen{} creates a tetrahedral mesh with a large number of slivers positioned on the faces of the regular cubic mesh (\cref{fig:kubik_smoothing:cube}).
After mesh optimization, all slivers are eliminated and the mesh quality is greatly improved (\cref{fig:kubik_smoothing:optimized,fig:kubik_smoothing:slivers,fig:kubik_smoothing:histo}).

The most difficult case is related to locked slivers, when a sliver chain has a number of fixed boundary faces and vertices. To solve this case, we suggest adding an additional prismatic layer between sliver chains and a fixed boundary, so that all sliver vertices are free to move.

\subsection{Testing the Anisotropic Encroachment Domain (AED) Based Algorithm}%
\label{sec:tests_encroachment}

The use of an anisotropic encroachment domain made it possible to reduce the excessive refinement of the surface meshes while keeping the mesh quality under control. For example, for the small AAU test structure (\cref{fig:aau_struct}), the usual DR from \tetgen{} produces a mesh with about 1 million elements in contrast to 40 thousand elements in the initial Ansys HFSS mesh.
 
\begin{figure}[p]%
  \subcaptionbox{\label{fig:dr_overrefinement:a}}{%
    \includegraphics[height = 0.32\textwidth]{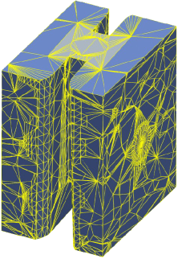}%
  }
  \hfill{}%
  \subcaptionbox{\label{fig:dr_overrefinement:b}}{%
    \includegraphics[height = 0.32\textwidth]{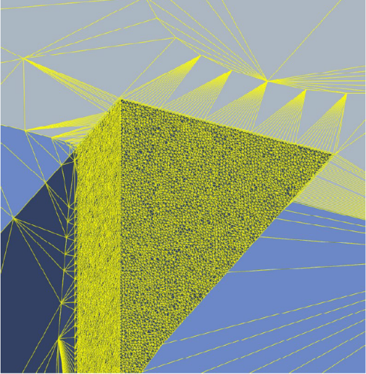}%
  }
  \hfill{}%
  \subcaptionbox{\label{fig:dr_overrefinement:c}}{%
    \includegraphics[height = 0.32\textwidth]{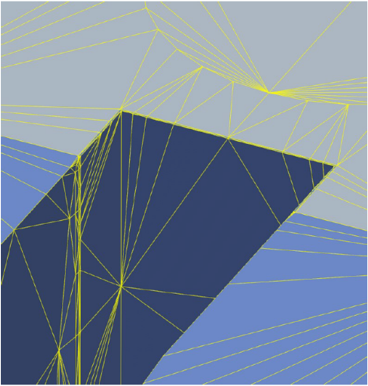}%
  }
  \caption{%
     Comparison of the classic DR and AED based refinement:
    \protect\subref{fig:dr_overrefinement:a} input structure,
    \protect\subref{fig:dr_overrefinement:b} the standard Delaunay refinement,
    and \protect\subref{fig:dr_overrefinement:c} AED based insertion rules}%
  \label{fig:dr_overrefinement}%
\end{figure}

The new AED based algorithm is demonstrated in \cref{fig:dr_overrefinement}.
\Cref{fig:dr_overrefinement:a} shows a very simple input structure. A very coarse initial mesh is refined using the DR algorithm in which an element is split either due to quality constraints of due to large size. \Cref{fig:dr_overrefinement:b} shows a fragment of the mesh near a very thin layer of material with an aspect ratio of about $100$: the mesh is unreasonably dense, the total number of mesh elements is about 1 million.
\Cref{fig:dr_overrefinement:c} shows a mesh obtained with our AED based approach: there is no overrefinement, the mesh contains about 10 thousand elements.

\begin{figure}[p]%
  \centering{}%
  \includegraphics[width = 0.88\textwidth]{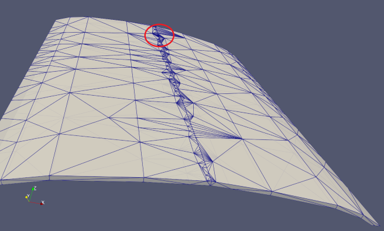}%
  \\[1ex]%
  \includegraphics[height = 0.393\textheight]{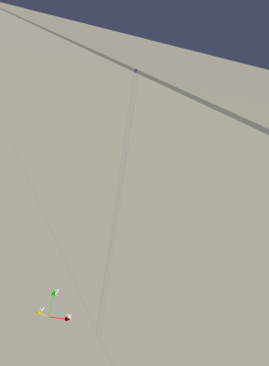}%
  \hspace{2ex}%
  \includegraphics[height = 0.393\textheight]{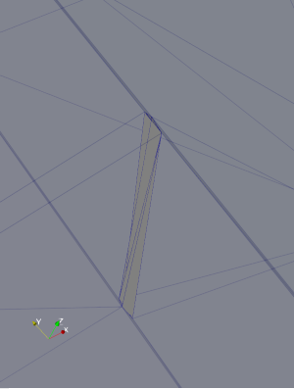}%
  \caption{A tiny step in the CAD model of the antenna array unit. Four tiny slivers with dihedral angles less than $10^{-2}$ are created during the AMR process}%
  \label{fig:aau_radome_step}%
\end{figure}

\begin{figure}[p]%
  \centering{}%
  \includegraphics[width=\textwidth]{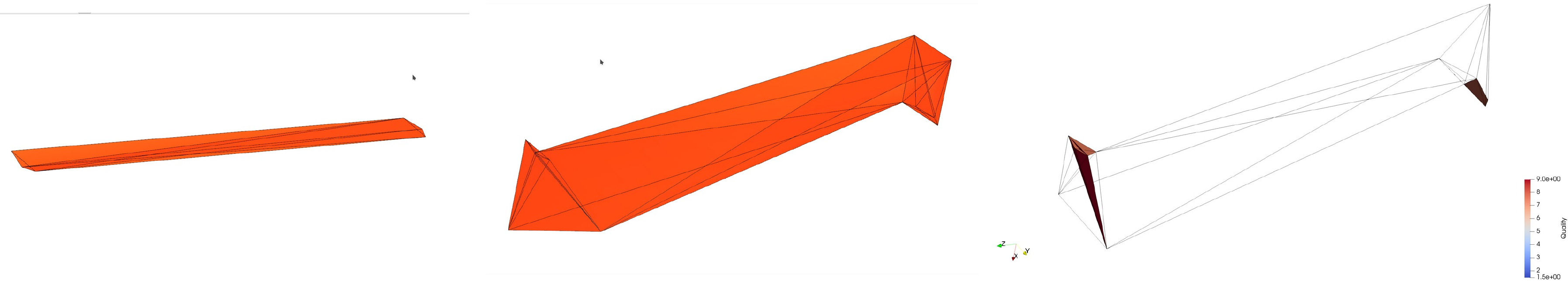}%
  \caption{Left: a locked chain of four slivers, small facets on the left and right are fixed. Center: tetrahedral mesh after adding near-boundary prisms is untangled and optimized. Right: the worst tetrahedra are coupled to the fixed boundary faces}%
  \label{fig:aau_locked_slivers}%
\end{figure}

In the AAU test case, there is a tiny step in the middle of a thin radome element, see \cref{fig:aau_radome_step}. Such defects often appear when using copy by reflection or copy by move in the geometry design stage. They are unavoidable, and therefore meshing algorithms used for serial simulations must handle such cases without manual healing of the geometry.
Thus, when using a DR-like insertion rule for mesh refinement, some post-processing is required to maintain mesh quality during the AMR process. This particular case is quite difficult because edge flips and mesh refinement are not able to remove this particular sliver chain.
To this end, we use a special procedure which pads the faces of the locked sliver by artificial prisms. Each boundary vertex of the sliver is padded by a newly added artificial edge which is surrounded by tetrahedra to keep the mesh consistent. This procedure turns locked slivers into unlocked ones, so that the optimization tool can be applied to eliminate them.
\Cref{fig:aau_locked_slivers} shows an enlarged and rotated view of the locked sliver chain from the AAU test case shown in \cref{fig:aau_radome_step}.
The success of the optimization procedure after adding an artificial prismatic layer near fixed faces is illustrated.

\section{Adaptive Mesh Refinement (AMR)}%
\label{sec:amr}

The general statement of the AMR problem is quite straightforward: given the current mesh and error indicators for each element (or vertex, or edge), and some limit on the number of elements, create a new mesh that provides less error in the FEM solution and has a number of elements not more than the prescribed limit.
In fact, a more detailed domain-specific requirements are usually specified for the AMR, for example:
\begin{itemize}
\item error equidistribution must be satisfied;
\item tetrahedra with largest errors must be split or element size in regions with largest errors must be reduced;
\item refinement must be computationally cheap compared to FEM solution;
\item most of the elements must remain unchanged;
\item refinement must be optimal in some sense.
\end{itemize}

There are several approaches to solve the AMR problem: bisection based algorithms, e.g.,~\cite{arnold2000locally}, Delaunay refinement~\cite{cheng2012delaunay}, and recovery of target mesh size distribution plus metric-based mesh adaptation by local mesh modifications of both topology and vertex positions~\cite{alliez2005variational, klingner2008aggressive, vasilevskii1999adaptive, agouzal1999adaptive}.

Experiments show that bisection is not suitable for industrial problems because it keeps badly shaped elements in the regions where they were present in the initial mesh. Also, bisection is not optimal because many points are inserted into the mesh just to maintain conformity but not because of large errors.

Metric-based optimization using both vertex smoothing and topological modifications seems to be the most powerful technique today, but it has a high computational cost.

Techniques similar to Delaunay refinement seem to be a good compromise for the AMR. The core part of the algorithm is simply point generation rules and point insertion into an existing tetrahedral mesh, which makes it cheap. At the same time, the original DR enjoys theoretical guarantees for mesh quality.

It is important to note that element-wise error indicators are not sufficient to fully control the AMR process. When a tetrahedron in the initial mesh is split, a whole neighborhood of this element is remeshed, and therefore one needs some criterion to decide which elements in the renewed mesh require splitting. The most natural way to do this is to somehow compute a target \emph{sizing function} (SF), which is the desired element diameter at each point. If the SF is set in vertices, then we can interpolate it for each new point inserted. Using this SF we can compute the target diameter of an element and decide whether we should split it or not.
Taking into account all the limitations and difficulties mentioned above, we have arrived at the general scheme of the AMR algorithm presented in \cref{alg:general,alg:tetsplit}.
In the algorithm,
\begin{description}
\item[$TargetSF(T')$] is computed by linear interpolation of the SF values at $T'$-s centroid (center of mass) from the values at other mesh vertices. For each new point $p$, $SF^{target}(p)$ is also computed by linear interpolation.

\item[$ratio = SF^{current} / SF^{target}$] is the ratio of the target and current sizing functions at some point or for some tetrahedron.

\item[$Queue$] is a priority queue, elements are sorted by ratio value, largest first.

\item[$TetSize(T)$] computes the 1D size of a tetrahedron, e.g., the maximum edge length or the edge length of a regular tetrahedron with the same volume as $T$.

\item[$SplitTet(T)$] generates a new point according to the point generation rules and inserts it into the mesh.
For the new tetrahedra from the refilled cavity, $ratio$ is computed and pairs $[T, ratio]$ are put in the priority queue.
\end{description}

\begin{algorithm}[t]%
  \caption{General AMR scheme}%
  \label{alg:general}%
  \begin{algorithmic}[1]
     \State $SF^{target} = ComputeTargetSF(Mesh,\,indicators)$ 
    \For{tet $T$ in $Mesh$} \Comment{Form the initial tet priority queue}
      \State $SF^{current}(T)$ = $TetSize(T)$ \Comment{Compute current size}
      \State $ratio = SF^{current}(T) / SF^{target}(T)$
      \State $Queue.Push([T,\,ratio])$
    \EndFor
    \While{$Queue$ is not empty and $Mesh.Size() < N_{max}$}
    \State [$T$,\,$ratio$] = $Queue.Pop()$
      \If{$T$ was not deleted before}
        \State $SplitTet(T)$
      \EndIf
    \EndWhile
    \State $OptimizeMesh()$ \Comment{Sliver removal and quality improvement}
  \end{algorithmic}%
\end{algorithm}%
\begin{algorithm}[t]%
  \caption{Procedure for tetrahedra splitting}%
  \label{alg:tetsplit}%
  \begin{algorithmic}[1]
    \Function{$SplitTet$}{tetrahedron $T$}
      $v = CircumCenter(T)$
      \If{$v$ encroaches upon face $f$}
        \State $v = CircumCenter(f)$
        \If{$v$ encroaches upon segment $s$}
          \State $v = Center(s)$
        \EndIf
      \EndIf
      \State $newTets = InsertPoint(v)$ \Comment{Delaunay-type insertion algorithm}
      \For{$T'$ in $newTets$} \Comment{Put new cavity tets in the priority queue}
        \State $SF^{current}(T') = TetSize(T')$
        \State $SF^{target}(T') = TargetSF(T')$
        \State $ratio = SF^{current}(T') / SF^{target}(T')$
        \State $Queue.Push([T',\,ratio])$
      \EndFor
    \EndFunction
  \end{algorithmic}%
\end{algorithm}%

The Delaunay-type insertion can still produce slivers by accident, so our mesh optimization tool is applied to eliminate them. Since the current version of the optimization tool does not use an adaptation metric, the problem of oversmoothing exists. To prevent oversmoothing, bad elements are put into the priority queue, and the worst elements with their small neighborhoods are optimized. This simple heuristics prevents excessive smearing of the size distribution in the adaptive mesh.

\section{Simulation Results}%
\label{sec:simulation_results}

Due to the complexity of the selected solutions, it is not possible to evaluate the accuracy of the numerical simulations by comparing the results with analytical solutions.
For this reason, in order to assess the efficiency of the proposed AMR pipeline, we tested it on two structures described above and compared the results with those obtained by Ansys HFSS 2020R1 industrial software. To make the experiments consistent, we used the same initial mesh generated by HFSS and the same heuristic convergence criteria: AMR stopps when the maximum difference in S-parameters between two subsequent solutions is less then the threshold $0.02$.
The limit for the increment of the number of elements at one refinement step was set to \qty{30}{\percent}.
It should be emphasized that the FEM solvers in our code and in HFSS provide the same solutions for the same mesh, which allows us to isolate the evaluation of AMR algorithms.

\Cref{tab:6cav,tab:aau} summarize the results: for both structures our approach gives satisfactory results, for the 6-cavity filter test we even get a significantly smaller mesh and a smaller number of AMR steps.
\begin{table}[p]%
   \begin{minipage}[b]{0.48\linewidth}
    \caption[short for lof]{AMR results for the resonator filter (6cav)}%
    \begin{tabular}{cccc}%
      \toprule%
      & initial size & final size & AMR steps\\%
      \midrule%
      HFSS      & \qty{9}{k} & \qty{391}{k} & 21\\%
      Our code  & \qty{9}{k} & \qty{330}{k} & 14\\%
      \bottomrule%
    \end{tabular}%
    \label{tab:6cav}%
   \end{minipage}
  \hfill{}%
   \begin{minipage}[b]{0.48\linewidth}
    \caption[short for lof]{AMR results for the AAU}%
    \label{tab:aau}%
    \begin{tabular}{cccc}%
      \toprule%
      & initial size & final size & AMR steps \\%
      \midrule%
      HFSS      & \qty{40}{k} & \qty{207}{k} & 7\\%
      Our code  & \qty{40}{k} & \qty{200}{k} & 7\\%
      \bottomrule%
    \end{tabular}%
  \end{minipage}
\end{table}

\begin{figure}%
  \centering{}%
  \includegraphics[width = 1.0\textwidth]{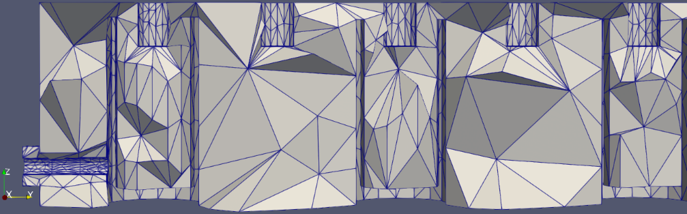}%
  \\[1ex]%
  \includegraphics[width = 1.0\textwidth]{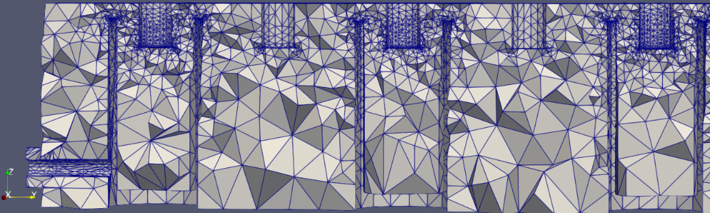}%
  \caption{Clips of the initial (top, \qty{9}{k} elements) and the final (bottom, \qty{330}{k} elements) meshes for the resonator filter 6cav}%
  \label{fig:6cav_ini_final_meshes}%
\end{figure}

\begin{figure}[p]%
  \centering{}%
  \includegraphics[width = 0.49\textwidth]{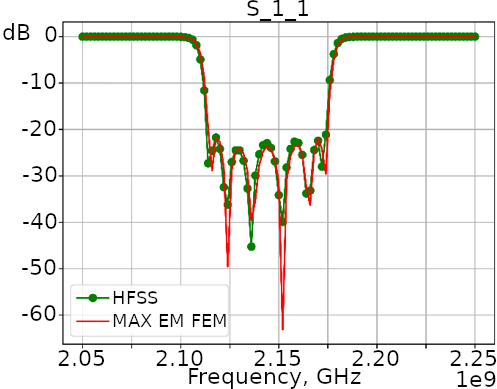}%
  \hfill{}%
  \includegraphics[width = 0.49\textwidth]{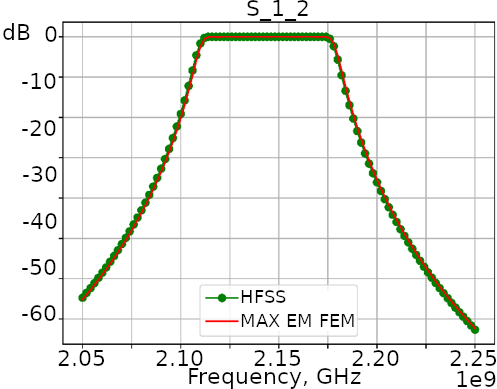}%
  \caption{Comparison of S-parameters computed on the final meshes: HFSS and current approach}%
  \label{fig:6cav_sparams}%
\end{figure}

\Cref{fig:6cav_ini_final_meshes} shows cuts of the initial and the final meshes for the 6cav structure. Note that the refinement here is extremely non-uniform, the mesh size is very small near edge singularities on some of the upper cylinders and quite large far from them. \Cref{fig:6cav_sparams} confirms that our approach perfectly reproduces the baseline solution.

\begin{figure}[p]%
  \centering{}%
  \includegraphics[height = 0.39\textwidth]{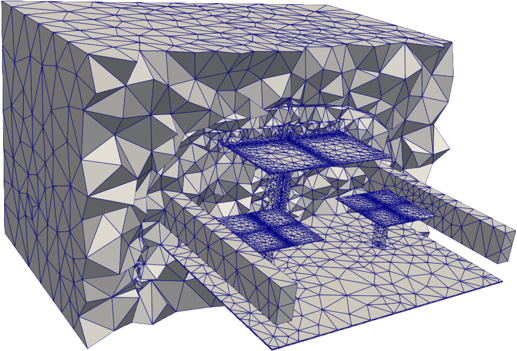}%
  \hfill{}%
  \includegraphics[height = 0.39\textwidth]{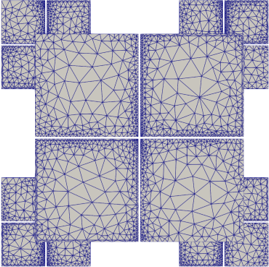}%
  \caption{Fragments of the final mesh for the AAU structure}%
  \label{fig:aau_meshes}%
\end{figure}

\begin{figure}[p]%
  \centering{}%
  \includegraphics[width = 0.49\textwidth]{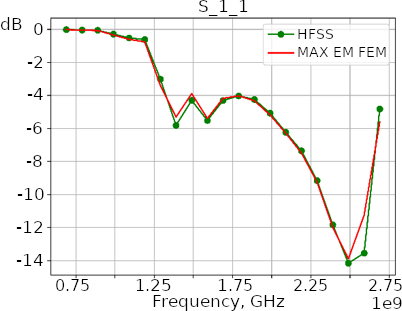}%
  \caption{S\textsubscript{11} parameter for the AAU case computed with adapted meshes}%
  \label{fig:aau_sparams}%
\end{figure}

Fragments of the final mesh for the AAU structure are shown in \cref{fig:aau_meshes}. Note that the surface mesh on square dielectric antenna lobes with thin metal stripes looks quasi-uniform. This is the merit of variational mesh smoothing.
\Cref{fig:aau_sparams} demonstrates good agreement of the S\textsubscript{11} parameter between our solver and HFSS.  Other parameters also show agreement.

\begin{figure}%
  \subcaptionbox{\label{fig:far_field:a}}{%
    \includegraphics[width=0.32\textwidth]{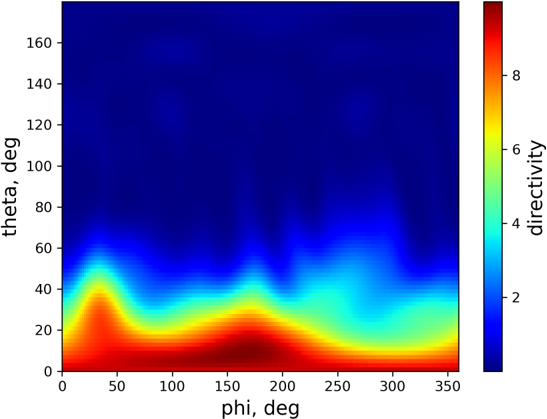}%
  }%
  \hfill{}%
  \subcaptionbox{\label{fig:far_field:b}}{%
    \includegraphics[width=0.32\textwidth]{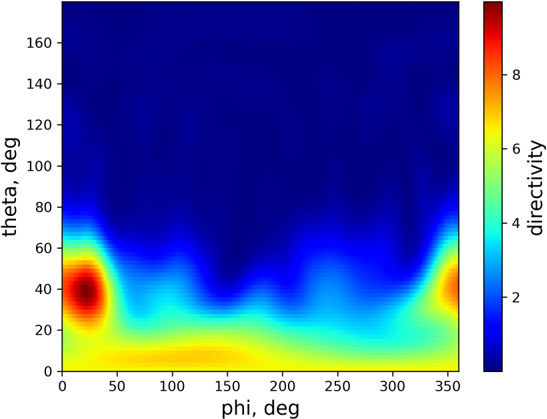}%
  }%
  \hfill{}%
  \subcaptionbox{\label{fig:far_field:c}}{%
    \includegraphics[width=0.32\textwidth]{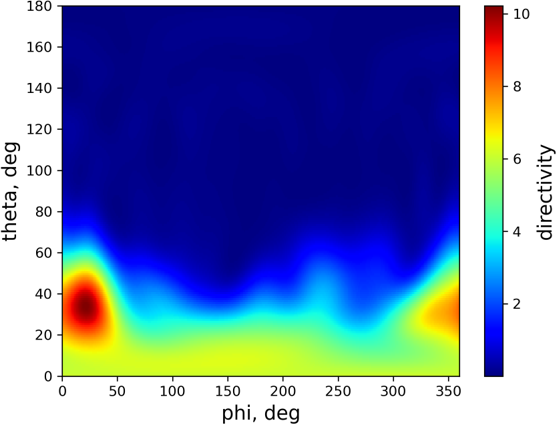}%
  }%
  \caption{Comparison of far-fields (directivity):
    \protect\subref{fig:far_field:a}~no smoothing,
    \protect\subref{fig:far_field:b}~with smoothing,
    \protect\subref{fig:far_field:c}~HFSS~AMR%
  }%
  \label{fig:far_field}%
\end{figure}

Apart from the S-parameters, the far-field intensity is also very important for antenna structures. We found that mesh smoothing strongly affects the accuracy (\cref{fig:far_field}).
At the same time, the increase in mesh quality is not so dramatic when comparing the quality distribution quality before and after smoothing.
This phenomenon should be further analyzed to understand how mesh quality affects solution accuracy.
Our current general assumption is that element quality can be particularly important in some regions, such as metallic sheets. Thus, a small increase in quality can significantly improve accuracy. Another observation is that S-parameters are less sensitive to mesh size and quality distribution than far-field intensity.

\section{Conclusion}%
\label{sec:conclusion}
An adaptive mesh refinement algorithmic toolchain for electromagnetic simulations was proposed and tested against the Ansys HFSS commercial baseline. Finite element error estimates based on the r.h.s.\ of the weak equations were used. The proposed approach uses only the order relation w.r.t.\ to error indicators to build the target mesh size distribution. The mesh refinement algorithm splits mesh elements ordered w.r.t.\ to the ratio of the current element size to the target size.
The splitting is done by point insertion using simple point generation rules and by extending a Delaunay-type insertion algorithm with anisotropic encroachment domain to prevent over-refinement inside thin material layers. Each AMR step is post-processed using a powerful variational smoothing algorithm. Computational results show that the developed toolchain can reach and even outperform the quality of commercial simulators.


\end{document}